\documentclass{article}

\usepackage{arxiv}

\usepackage[utf8]{inputenc} 
\usepackage[T1]{fontenc}    
\usepackage{hyperref}       
\usepackage{url}            
\usepackage{booktabs}       
\usepackage{amsfonts}       
\usepackage{nicefrac}       
\usepackage{microtype}      
\usepackage{subcaption}
\usepackage{mathrsfs}
\usepackage[framed,numbered,autolinebreaks,useliterate]{mcode} 
\usepackage{comment}
\usepackage{lipsum}
\usepackage{graphicx}
\usepackage{amsmath}
\usepackage{color}
\usepackage{colortbl}
\usepackage{xcolor}
\usepackage{transparent}
\graphicspath{{Bilder/}}
\usepackage{tikz}
\usetikzlibrary{arrows}

\DeclareMathOperator*{\argminA}{arg\,min}

\title{Parameter estimation of temperature dependent material parameters in the cooling process of TMCP steel plates}

\author{
  Dimitri Rothermel \\
  Department of Numerical Mathematics\\
  Saarland University\\
  Saarbr\"ucken, Germany \\
  \texttt{dimitri.rothermel@num.uni-sb.de} \\
       \And
  Thomas Schuster \\
  Department of Numerical Mathematics\\
  Saarland University\\
  Saarbr\"ucken, Germany \\
  \texttt{thomas.schuster@num.uni-sb.de} \\
     \And
  Roland Schorr \\
  Research \& Development\\
  AG der Dillinger H\"uttenwerke\\
  Dillingen/Saar, Germany \\
  \texttt{roland.schorr@dillinger.biz} \\
     \And
 Martin Peglow \\
  Research \& Development\\
  AG der Dillinger H\"uttenwerke\\
  Dillingen/Saar, Germany \\
  \texttt{martin.peglow@dillinger.biz} \\
}

\begin{document}
\maketitle

\begin{abstract}
Accelerated cooling is a key technology in producing thermomechanically controlled processed (TMCP) steel plates. In a TMCP process hot plates are subjected to a strong cooling what results in a complex microstructure leading to increased strength and fracture toughness. The microstructure is strongly affected by the temperature evolution during the cooling process as well as residual stresses and flatness deformations. Therefore, the full control (quantification) of the temperature evolution is very important regarding plate design and processing. It can only be achieved by a thermophysical characterization of the material and the cooling system.  In this paper, we focus on the thermophysical characterization of the material parameters.
Mathematically, we consider a specific inverse heat conduction problem. The temperature evolution of a heated steel plate passing through the cooling device is modeled by a 1D nonlinear partial differential equation (PDE) with unknown temperature dependent material parameters, which describe the characteristics of the underlying material.
We present a numerical approach to identify these material parameters up to some canonical ambiguity without any a priori information. 
\end{abstract}

\newpage

\section{Introduction}\label{sec:eins}
This paper is motivated by the cooling process in the production of heavy plates made of steel.
From a metallurgical point of view it is known that the specific temperature profile in a cooling device can affect the crystalline microstructure of steel and is therefore crucial to meet various customer requirements.
That is why understanding and modeling the cooling process of heated steel plates is of great interest in adjusting the mechanical properties (toughness, strength, weldability, etc.) of the final product.
To achieve a full control of the cooling process however, we first need to be able to thermophysically characterize the underlying material. The application of the cooling water on the hot surfaces of the heavy plates leads to a heat extraction only on the boundary, whereas the temperature evolution \textit{inside} of the plate is solely driven by its material parameters, i.e. the volumetric heat capacity and the thermal conductivity, which we denote by $C$ and $k$, respectively. 

There is already a vast number of publications focusing on the determination of those material parameters that have been presented in similar contexts, involving different measurement data setups to formulate so-called inverse heat conduction problems (IHCP). 

In \cite{egger} for example, the authors proved the uniqueness of solutions concerning a model with constant Neumann boundary conditions and a single boundary temperature measurement by extending Cannon's uniqueness result (see \cite{cannon}) for the non-stationary case. Rather theoretical cases, e.g. for homogenuous inital temperatures and $k$ being a multiple of $C$, were also discussed by Cannon et al., see \cite{cannonDuchateau}.
Further research results concerning the existence and uniqueness of solutions, as well as some numerical analyis concerning the quenching process for a different model with homogenuous Dirichlet boundary conditions can be found in \cite{rincon1}, \cite{rincon2} and the references contained therein. 

The extensive literature also contains various models where the material parameters only depend on the space variable $z$ (see \cite{benyuZou}, \cite{engl}) or only on the time variable $t$, see \cite{hussein}. The estimation of temporally \textit{and} spatially varying diffusion coefficients is considered in \cite{kunisch}, see also the references contained therein.

In the 1D heat conduction PDE model considered in this article
\begin{align} \label{introductionModel}
C(u)u_t = (k(u)u_z)_z,
\end{align}
which we discuss in Section \ref{sec:drei} in more detail, the material parameters depend on the temperature $u(t,z)$ itself due to the high temperature gradients and the phase changes in the material. Note, that \ref{introductionModel} is often also used in the modeling of heating processes, i.e. the addition of heat, e.g. in simulations of blast furnaces or in space research. However, in our case we are interested in the cooling process, i.e. the extraction of heat from heated heavy plates made of steel by applying cooling water on the surfaces.

Numerous methods for solutions to IHCPs with similar models containing temperature dependent material parameters were also discussed and analyzed. In \cite{hankeScherzer} for example, the authors used a conjugate gradient based equation error method by solving a system of linear operator equations for the diffusion parameter. Cui et al. focused on the sensitivity analysis of their underlying gradient based method by using the complex-variable differentiation method (CVDM), see \cite{cuiGaoZhang}. There, the authors present simple numerical results, showing that they are able to identify the temperature dependent material parameters, i.e. the coefficients of second order polynomials. In \cite{huangYan}, the authors also try to determine more complex material parameters $k(u)$ and $C(u)$ simultaneously by arguing, that once the temperature distribution $u(t,z)$ is known one can replace $k(u)$ by $k(t,z)$ and $C(u)$ by $C(t,z)$. Yet, using a numerical scheme, the determination of some time and space dependent material parameters gets computationally more expensive when refining the spatial and temporal grid size, which obviously is a drawback. In this paper, we want to propose a method to solve simultaneously for general functions $k(u)$ and $C(u)$ in the space $\mathcal{C}^1$ of continuously differentiable functions without any a priori information and functional forms that are rather general.

The focus of this paper is application-oriented. After introducing the experimental setup in Section \ref{sec:zwei} and the mathematical description in Section \ref{sec:drei}\&\ref{sec:vier}, we want to show a practical guide in Section \ref{sec:fuenf} on how to numerically implement the underlying parameter estimation method.


\section{The experimental setup}\label{sec:zwei}
\begin{figure}[h]
    \centering
        \def\svgwidth{450pt}
\begingroup%
  \makeatletter%
  \providecommand\color[2][]{%
    \errmessage{(Inkscape) Color is used for the text in Inkscape, but the package 'color.sty' is not loaded}%
    \renewcommand\color[2][]{}%
  }%
  \providecommand\transparent[1]{%
    \errmessage{(Inkscape) Transparency is used (non-zero) for the text in Inkscape, but the package 'transparent.sty' is not loaded}%
    \renewcommand\transparent[1]{}%
  }%
  \providecommand\rotatebox[2]{#2}%
  \newcommand*\fsize{\dimexpr\f@size pt\relax}%
  \newcommand*\lineheight[1]{\fontsize{\fsize}{#1\fsize}\selectfont}%
  \ifx\svgwidth\undefined%
    \setlength{\unitlength}{1998.57202148bp}%
    \ifx\svgscale\undefined%
      \relax%
    \else%
      \setlength{\unitlength}{\unitlength * \real{\svgscale}}%
    \fi%
  \else%
    \setlength{\unitlength}{\svgwidth}%
  \fi%
  \global\let\svgwidth\undefined%
  \global\let\svgscale\undefined%
  \makeatother%
  \begin{picture}(1,0.30475087)%
    \lineheight{1}%
    \setlength\tabcolsep{0pt}%
    \put(0,0){\includegraphics[width=\unitlength,page=1]{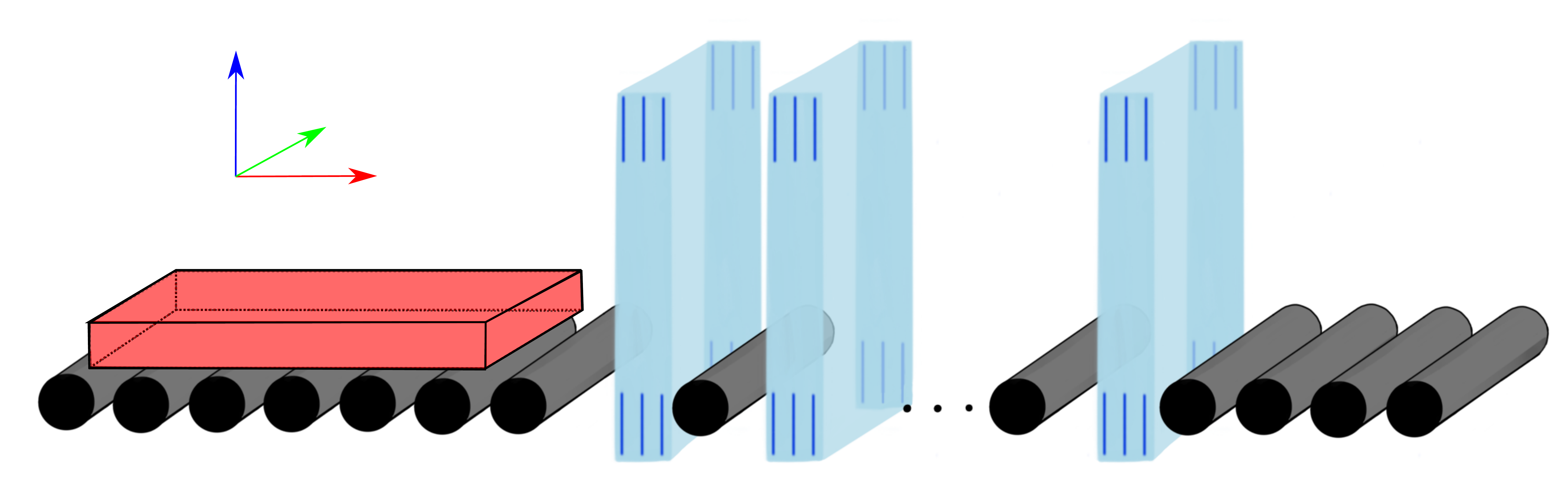}}%
    \put(0.24447497,0.18634247){\color[rgb]{1,0,0}\makebox(0,0)[lt]{\lineheight{1.25}\smash{\begin{tabular}[t]{l}$x$\end{tabular}}}}%
    \put(0.20931058,0.22867846){\color[rgb]{0,1,0}\makebox(0,0)[lt]{\lineheight{1.25}\smash{\begin{tabular}[t]{l}$y$\end{tabular}}}}%
    \put(0.14166512,0.27837988){\color[rgb]{0,0,1}\makebox(0,0)[lt]{\lineheight{1.25}\smash{\begin{tabular}[t]{l}$z$\end{tabular}}}}%
  \end{picture}%
\endgroup%

    \captionof{figure}{Sketch of the laminar cooling system}
    \label{fig:kuehlanlage}
\end{figure}
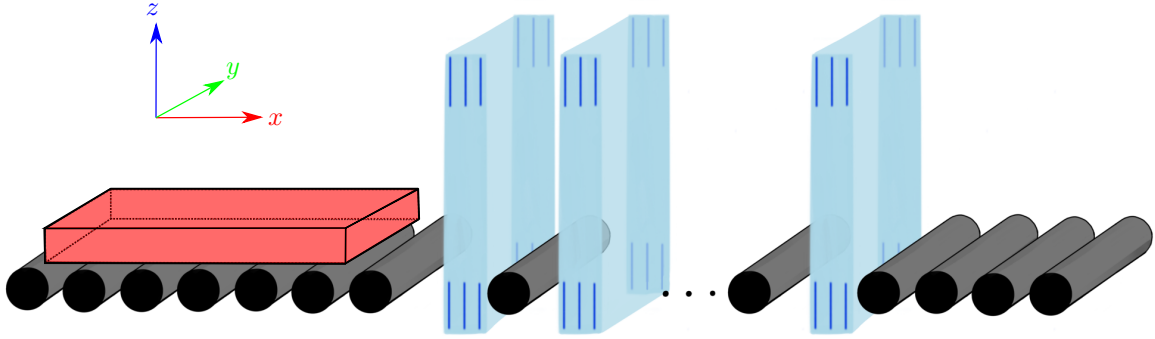

In order to describe the underlying laminar cooling system as necessary as possible, we refer to the sketch in Figure \ref{fig:kuehlanlage}. Multiple conveyor rollers transport the heated heavy plate in $x$-direction towards the cooling system, consisting of several water cooling zones that are indicated by dots. 

Depending on the material parameters (thermal conductivity, volumetric heat capacity), applying cooling water on top and bottom surface (with respect to the $x$-$y$-plane) leads to a special heat transfer within the plate and thus to a very specific final temperature profile.

\begin{figure}[h]
    \centering
      \def\svgwidth{270pt}
\begingroup%
  \makeatletter%
  \providecommand\color[2][]{%
    \errmessage{(Inkscape) Color is used for the text in Inkscape, but the package 'color.sty' is not loaded}%
    \renewcommand\color[2][]{}%
  }%
  \providecommand\transparent[1]{%
    \errmessage{(Inkscape) Transparency is used (non-zero) for the text in Inkscape, but the package 'transparent.sty' is not loaded}%
    \renewcommand\transparent[1]{}%
  }%
  \providecommand\rotatebox[2]{#2}%
  \newcommand*\fsize{\dimexpr\f@size pt\relax}%
  \newcommand*\lineheight[1]{\fontsize{\fsize}{#1\fsize}\selectfont}%
  \ifx\svgwidth\undefined%
    \setlength{\unitlength}{553.75032065bp}%
    \ifx\svgscale\undefined%
      \relax%
    \else%
      \setlength{\unitlength}{\unitlength * \real{\svgscale}}%
    \fi%
  \else%
    \setlength{\unitlength}{\svgwidth}%
  \fi%
  \global\let\svgwidth\undefined%
  \global\let\svgscale\undefined%
  \makeatother%
  \begin{picture}(1,0.7160506)%
    \lineheight{1}%
    \setlength\tabcolsep{0pt}%
    \put(0,0){\includegraphics[width=\unitlength,page=1]{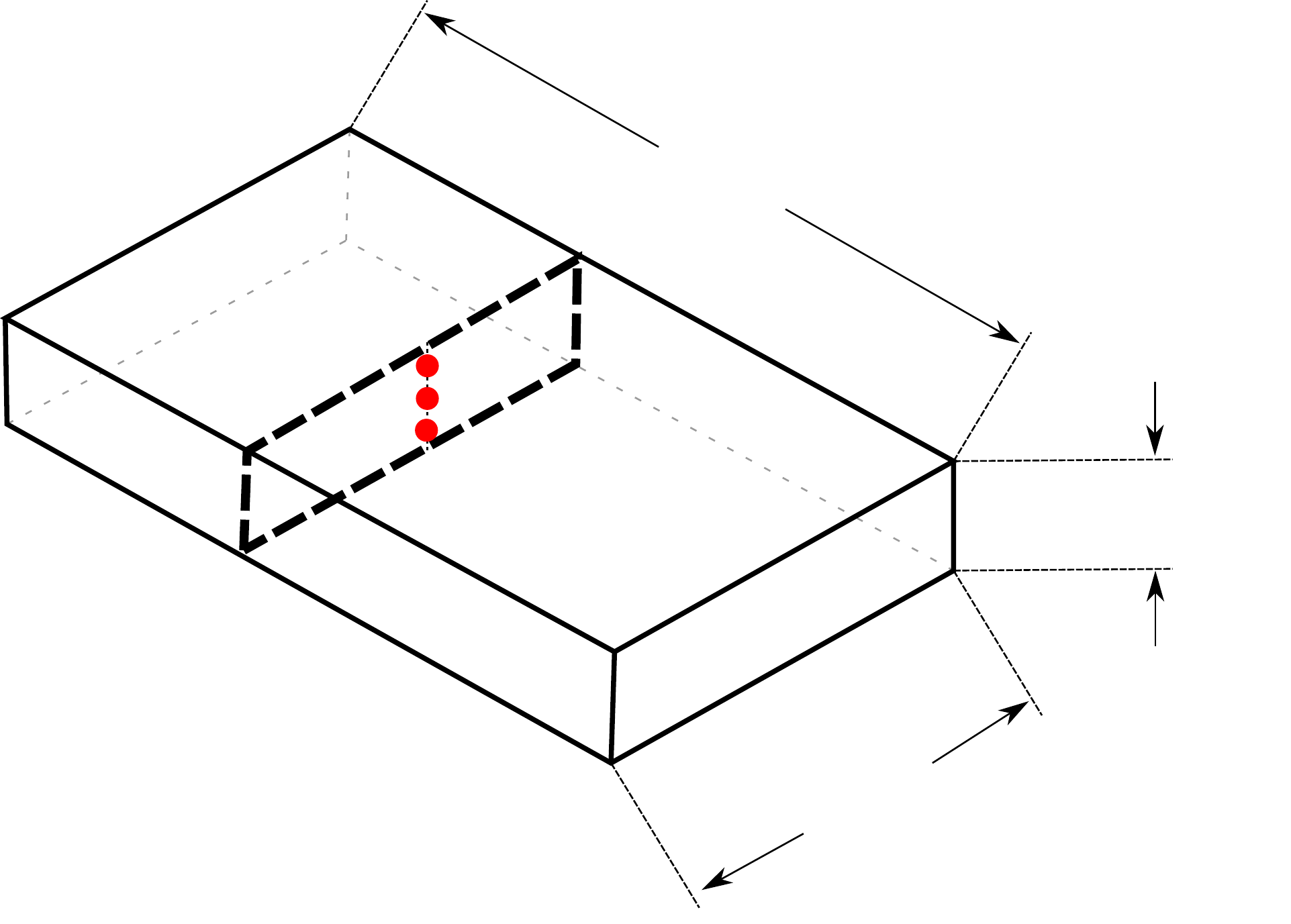}}%
    \put(0.50842584,0.56520439){\color[rgb]{0,0,0}\makebox(0,0)[lt]{\lineheight{1.25}\smash{\begin{tabular}[t]{l}Length\end{tabular}}}}%
    \put(0.75381997,0.30115106){\color[rgb]{0,0,0}\makebox(0,0)[lt]{\lineheight{1.25}\smash{\begin{tabular}[t]{l}Thickness\end{tabular}}}}%
    \put(0.62636169,0.08224681){\color[rgb]{0,0,0}\makebox(0,0)[lt]{\lineheight{1.25}\smash{\begin{tabular}[t]{l}Width\end{tabular}}}}%
    \put(0,0){\includegraphics[width=\unitlength,page=2]{plate1.pdf}}%
    \put(0.18139486,0.00954959){\color[rgb]{1,0,0}\makebox(0,0)[lt]{\lineheight{1.25}\smash{\begin{tabular}[t]{l}$x$\end{tabular}}}}%
    \put(0.17547788,0.20241707){\color[rgb]{0,1,0}\makebox(0,0)[lt]{\lineheight{1.25}\smash{\begin{tabular}[t]{l}$y$\end{tabular}}}}%
    \put(0.01228168,0.29841473){\color[rgb]{0,0,1}\makebox(0,0)[lt]{\lineheight{1.25}\smash{\begin{tabular}[t]{l}$z$\end{tabular}}}}%
  \end{picture}%
\endgroup%

    \captionof{figure}{Placement of thermocouples in the heavy plate}
    \label{fig:platte_duenn}
    \end{figure}

\newpage 
Typically, the heavy plates under consideration have a very low thickness in comparison to their length and width, see Figure \ref{fig:platte_duenn} for a skewed representation of the heavy plate dimensions. Usually, the length is $10$ to $30$m, whereas the thickness is around $30$ to $200$mm.

Therefore, the main information about the temperature evolution while cooling lies in the heat transfer in z-direction. 
Hence, we need to measure temperatures at different depths with respect to the thickness of the heavy plate, i.e. at locations with different $z$-components.
We achieve this by placing three so-called thermocouples into the heavy plate, indicated by the red dots in Figure \ref{fig:platte_duenn}. 
Note, that in our experimental setup the thermocouples are positioned right in the middle with respect to width and also symmetrically with respect to thickness, i.e. we have one core thermocouple and two near-surface thermocouples, one at the top and one at the bottom. 

\begin{center}
\begin{minipage}{0.5\linewidth}
    \def\svgwidth{350pt}
\begingroup%
  \makeatletter%
  \providecommand\color[2][]{%
    \errmessage{(Inkscape) Color is used for the text in Inkscape, but the package 'color.sty' is not loaded}%
    \renewcommand\color[2][]{}%
  }%
  \providecommand\transparent[1]{%
    \errmessage{(Inkscape) Transparency is used (non-zero) for the text in Inkscape, but the package 'transparent.sty' is not loaded}%
    \renewcommand\transparent[1]{}%
  }%
  \providecommand\rotatebox[2]{#2}%
  \newcommand*\fsize{\dimexpr\f@size pt\relax}%
  \newcommand*\lineheight[1]{\fontsize{\fsize}{#1\fsize}\selectfont}%
  \ifx\svgwidth\undefined%
    \setlength{\unitlength}{772.35864535bp}%
    \ifx\svgscale\undefined%
      \relax%
    \else%
      \setlength{\unitlength}{\unitlength * \real{\svgscale}}%
    \fi%
  \else%
    \setlength{\unitlength}{\svgwidth}%
  \fi%
  \global\let\svgwidth\undefined%
  \global\let\svgscale\undefined%
  \makeatother%
  \begin{picture}(1,0.32371112)%
    \lineheight{1}%
    \setlength\tabcolsep{0pt}%
    \put(0,0){\includegraphics[width=\unitlength,page=1]{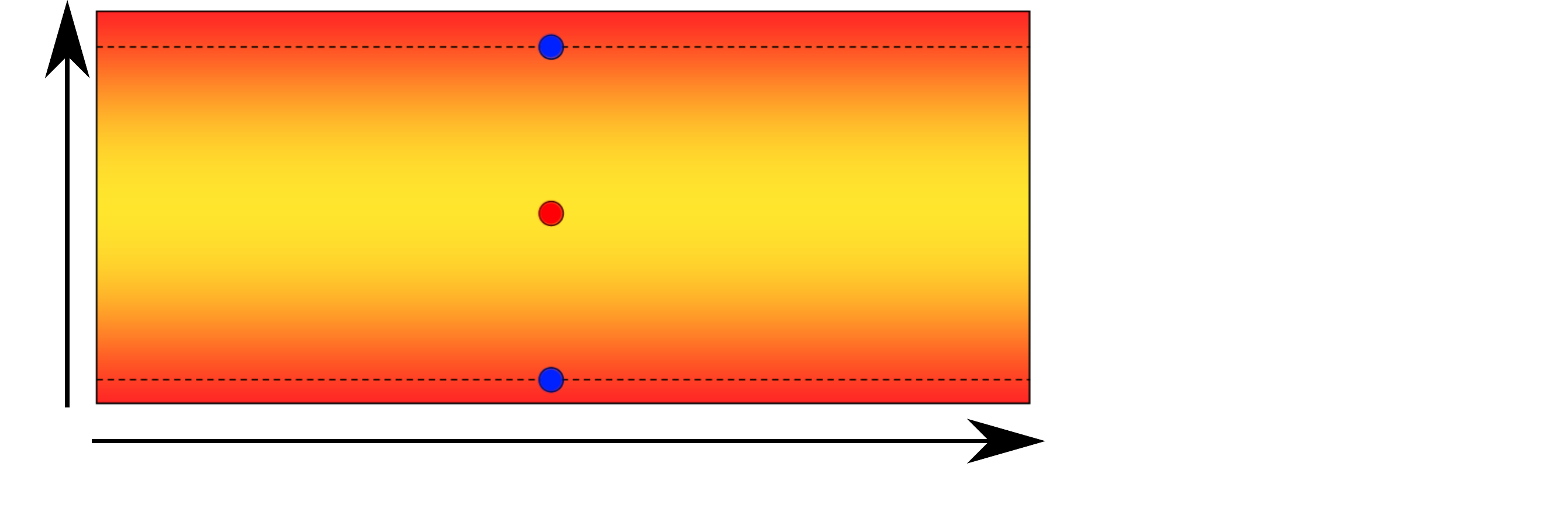}}%
    \put(0.25061471,0.02763866){\color[rgb]{0,0,0}\makebox(0,0)[lt]{\begin{minipage}{0.20253357\unitlength}\raggedright \end{minipage}}}%
    \put(0.317756,0.00058794){\color[rgb]{0,0,0}\makebox(0,0)[lt]{\lineheight{1.25}\smash{\begin{tabular}[t]{l}Width\end{tabular}}}}%
    \put(0.02951086,0.10087121){\color[rgb]{0,0,0}\rotatebox{90}{\makebox(0,0)[lt]{\lineheight{1.25}\smash{\begin{tabular}[t]{l}Thickness\end{tabular}}}}}%
    \put(0,0){\includegraphics[width=\unitlength,page=2]{plate2.pdf}}%
    \put(0.66816679,0.07896571){\color[rgb]{0,0,0}\makebox(0,0)[lt]{\lineheight{1.25}\smash{\begin{tabular}[t]{l}$z=0$\end{tabular}}}}%
    \put(0,0){\includegraphics[width=\unitlength,page=3]{plate2.pdf}}%
    \put(0.66521484,0.17948464){\color[rgb]{0,0,0}\makebox(0,0)[lt]{\lineheight{1.25}\smash{\begin{tabular}[t]{l}$z=\frac{L}{2}$\end{tabular}}}}%
    \put(0.66454982,0.28149931){\color[rgb]{0,0,0}\makebox(0,0)[lt]{\lineheight{1.25}\smash{\begin{tabular}[t]{l}$z=L$\end{tabular}}}}%
    \put(0,0){\includegraphics[width=\unitlength,page=4]{plate2.pdf}}%
  \end{picture}%
\endgroup%

    \captionof{figure}{Color-coded thermocouples at different depths $z$}
    \label{fig:ausschnitt}
\end{minipage}%
\end{center}

Figure \ref{fig:ausschnitt} shows a larger section of the marked cross-section in Figure \ref{fig:platte_duenn}. This time, the thermocouples are color-coded in such a way, that we can directly assign the measured temperature curves in Figure \ref{fig:kurven} to its corresponding depth by color. The variable $L>0$ is the length of the distance between the top and bottom thermocouple. By applying cooling water evenly on top ($z>L$) and bottom surface ($z<0$) of the heavy plate, we can assume that for some fixed depth $0<z<L$ the temperatures are equal with respect to width at any time.

\begin{center}
\def\svgwidth{350pt}
\begingroup%
  \makeatletter%
  \providecommand\color[2][]{%
    \errmessage{(Inkscape) Color is used for the text in Inkscape, but the package 'color.sty' is not loaded}%
    \renewcommand\color[2][]{}%
  }%
  \providecommand\transparent[1]{%
    \errmessage{(Inkscape) Transparency is used (non-zero) for the text in Inkscape, but the package 'transparent.sty' is not loaded}%
    \renewcommand\transparent[1]{}%
  }%
  \providecommand\rotatebox[2]{#2}%
  \newcommand*\fsize{\dimexpr\f@size pt\relax}%
  \newcommand*\lineheight[1]{\fontsize{\fsize}{#1\fsize}\selectfont}%
  \ifx\svgwidth\undefined%
    \setlength{\unitlength}{668.06421bp}%
    \ifx\svgscale\undefined%
      \relax%
    \else%
      \setlength{\unitlength}{\unitlength * \real{\svgscale}}%
    \fi%
  \else%
    \setlength{\unitlength}{\svgwidth}%
  \fi%
  \global\let\svgwidth\undefined%
  \global\let\svgscale\undefined%
  \makeatother%
  \begin{picture}(1,0.45683233)%
    \lineheight{1}%
    \setlength\tabcolsep{0pt}%
    \put(0,0){\includegraphics[width=\unitlength,page=1]{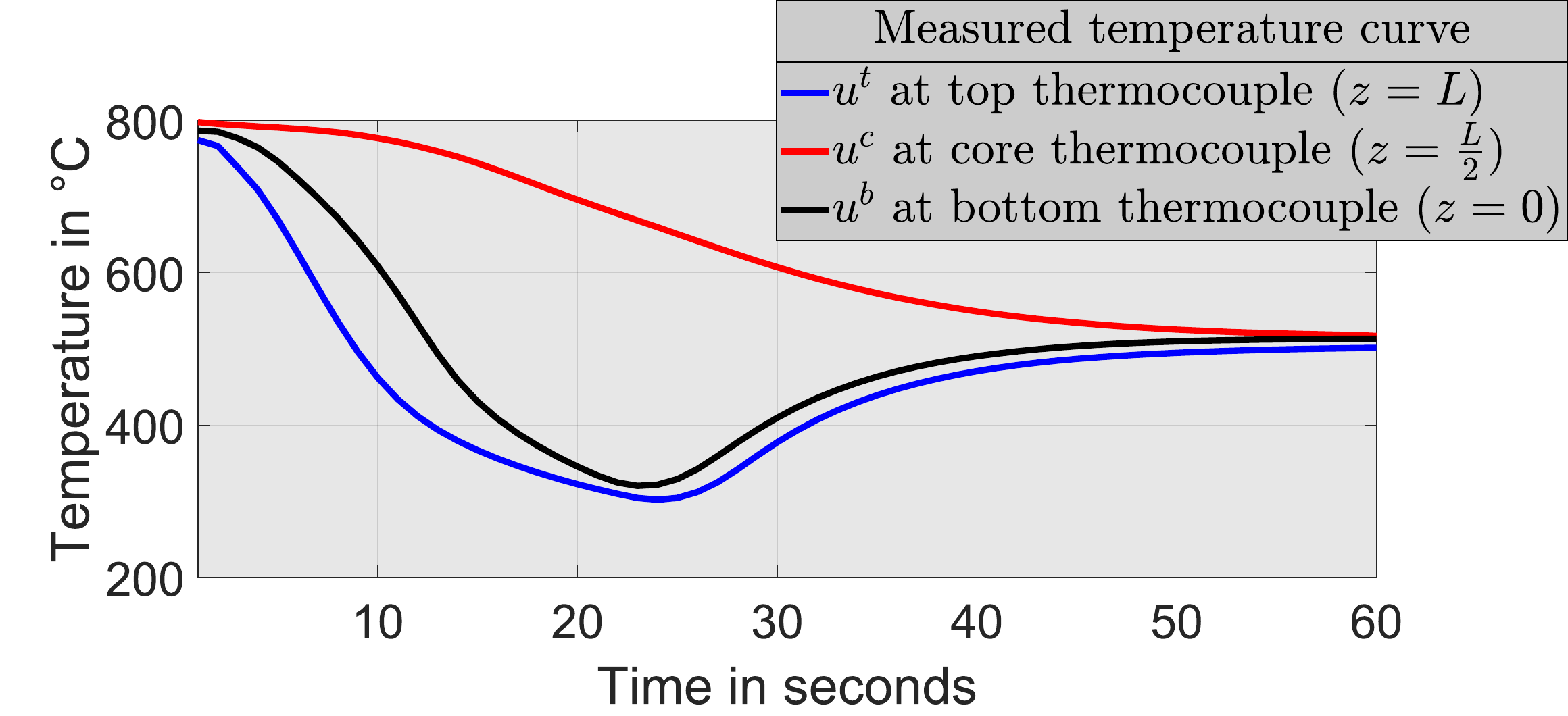}}%
  \end{picture}%
\endgroup%

\captionof{figure}{Measured temperatures at different thermocouple locations}
\label{fig:kurven}
\end{center}

In Figure \ref{fig:kurven}, the measured temperatures $u^t, u^c$ and $u^b \in \mathbb{R}_+^m\  (m>0)$ are plotted as interpolated curves. They encode the history of a typical cooling scenario of such plates as follows:
\begin{itemize}

\item Before the cooling, at time $t=0$ [$s$], the heavy plate is almost evenly heated at approx. $780\  ^{\circ}\mathcal{C}$. 
\item During the water cooling process
($0<t<25\  [s]$) the temperature decreases whereas the cooling rate of the near-surface regions is larger than the cooling rate of the core region due to the smaller distance to the cooled surfaces. 
\item When leaving the cooling system the cooling process is getting interrupted and recovery starts, at times t > 25 [s]. Due to heat conduction under air cooling conditions heat flows from the core region to the near-surface regions, which can be seen from an increase of temperature (see $u^t$ and $u^b$) .
\end{itemize}


The exact dependence of the temperature on space and time is highly affected by the material parameters. In the following, these material parameters will be determined as a solution of an inverse heat conduction problem, where $u^t$ and $u^b$ are considered as Dirichlet boundary conditions of our PDE model and $u^c$ as the measured data, see Section 4 and 5.



\section{Mathematical Preliminaries} \label{sec:drei}
Before we can formulate the inverse heat conduction problem in the form of the parameter estimation of the temperature dependent material parameters, we must clarify which mathematical model we use. Also, we need to define the so-called \textit{forward operator} and the \textit{observation operator}.
    \subsection{Model description} 

As already mentioned in Section \ref{sec:zwei}, the heavy plates under consideration have a very small thickness compared to their length and width. Cooling of top and bottom surface of such plates results in a 1D heat conduction problem. The temperature evolution of every 2 points of this 1D object includes the thermophysical information of the regime in between these 2 points.

    For $L>0$, let $\Omega:=(0, L)$ contain the spatial variables $z$ which describe the positions/depths in the heavy plate in direction of thickness, while $z=0$ corresponds to the position of the bottom and $z=L$ to the position of the top thermocouple, see also Figure \ref{fig:ausschnitt}.
    
    Our 1D model for heat conduction in the time frame $I:=(0,T]$ for some $T>0$ is given then by the initial boundary value problem
    \begin{align}\label{pde1}
        \tilde{C}(u)u_t &= (\tilde{k}(u)u_z)_z, &(t,z)\in I \times \Omega \\
        u&= u^b, &t\in I, z=0, \label{pde2}\\\label{pde3}
        u&= u^t, &t\in I, z=L,\\ \label{pde4}
        u&= u_0, &t=0, z\in \Omega,
    \end{align}
    where the unknown solution $u:\bar{I}\times \bar{\Omega} \to \mathbb{R}_+$ represents the temperature at times $t\in[0,T]$ and depths $z \in [0,L]$. The function $u_0(z)$ is the known initial temperature distribution before cooling, i.e. for $t=0$. For the Dirichlet boundary conditions \eqref{pde2}-\eqref{pde3} we use the measured boundary temperatures $u^b$ and $u^t$, see Figure \ref{fig:kurven}. For brevity, we often omit the variables $t$, $z$ or $(t,z)$ and, e.g., write $u$ instead of $u(t,z)$. By the subscripts $t$ and $z$ we refer to the derivatives $\frac{d}{dt}$ and $\frac{d}{dz}$ with respect to time and space, respectively.

    Let $U=[u_{min},u_{max}]$ with $0\leq u_{min}<u_{max}<\infty$ be the interval that covers all occurring temperatures $u$ in \eqref{pde1}-\eqref{pde4}. In our case, i.e. the cooling of heated heavy plates, the temperature of the plates are maximal at the start of the cooling process and do not fall below the temperature of the cooling water. Thus, we can set, e..g.,
    \begin{align*}
    u_{min} = 0 \text{\ \ \    and \ \ \   } u_{max} = \max\limits_{z \in \bar{\Omega}}\  u_0(z).
\end{align*}      
The functions 
    \begin{align*}
    \tilde{C}:& U\subset \mathbb{R}_+ \to \mathbb{R}_+\\\tilde{k}:& U\subset \mathbb{R}_+ \to \mathbb{R}_+    
    \end{align*}
    denote the {\em volumetric heat capacity} and the {\em thermal conductivity}, respectively. Let us assume that $\tilde{k},\tilde{C}\in \mathcal{C}^1(U)$. These functions represent the material parameters which lead to a unique heat conduction behaviour and consequently to a very specific temperature solution $u(t,z),\ (t,z)\in \bar{I} \times \bar{\Omega}$. Both functions depend on $u \in \mathbb{R}_+$ itself due to the presence of high temperature gradients while cooling and the phase transitions in the crystalline microstructure of the material. Thus, the 1D heat equation \eqref{pde1} is nonlinear in $u$.

    For this paper, we aim to put the focus on the parameter estimation methodology consisting of the determination of these material parameters $\tilde{C}$ and $\tilde{k}$ from core temperature knowledge rather than discussing the mathematical aspects of abstract function spaces and the theory contained therein. For that reason, we want to make the generalizing assumption that all functions in \eqref{pde1}-\eqref{pde4} are sufficiently smooth in the sense that the heat equation model is well-defined. Furthermore, we would also like to assume that a {\em unique} solution $\mathcal{U}\ni u:\bar{I}\times \bar{\Omega} \to \mathbb{R}_+$ always {\em exists}. Without specifying the function space $\mathcal{U}$, we note that the point evaluation is valid in time and space, i.e. $u(t,z)$ is at least continuous in $t \in \bar{I}$ and $z \in \bar{\Omega}$. We recommend \cite{roubicek} for more information about the abstract function spaces.
    
\subsection{The forward operator}
Let $\mathscr{K}$ and $\mathscr{C}$ be the spaces of all admissible functions $\tilde{k}$ and $\tilde{C}$ in \eqref{pde1}, respectively. 
We first define the operator
\begin{align}
    \label{wrongforwardOperator1}
    \tilde{F}:\mathscr{K}\times \mathscr{C}&\to \mathscr{U},\\   \label{wrongforwardOperator2}
     (\tilde{k},\tilde{C}) &\mapsto u,
\end{align}
where $u\in \mathscr{U}$ is the solution to the initial boundary value problem \eqref{pde1}-\eqref{pde4} corresponding to material parameters $(\tilde{k},\tilde{C})$. We face problems in defining some forward operator as in \eqref{wrongforwardOperator1}-\eqref{wrongforwardOperator2} explicitly, because the material parameters $\tilde{k}$ and $\tilde{C}$ depend on the temperature $u$, which in turn depends on $\tilde{k}$ and $\tilde{C}$. It is possible to formulate an implicitly defined forward operator, but that is not furthermore pursued in this article. Rather, we want to get rid of the temperature dependency of the input functions without discarding the underlying PDE, i.e. we want to retain the heat conduction model \eqref{pde1}-\eqref{pde4} up to some slight changes.

For that reason, let $\pi_n:u_{min}=u_1<u_2<\cdots<u_n=u_{max}$ be a partition of the interval $U=[u_{min},u_{max}]$. Moreover, let $\underline{k}=(k_1,\dots,k_n)^T\in\mathbb{R}^n_+$ and $\underline{C}=(C_1,\dots,C_n)^T\in\mathbb{R}^n_+$ be two sets of values. Given $\pi_n$, $\underline{k}$ and $\underline{C}$ we can construct piecewise cubic functions $k$ and $C$ in $\mathcal{C}^1(U)$, such that 
\begin{align*}
        k(u_i)=k_i,\\
    C(u_i)=C_i,
\end{align*}
for $i=1,\dots,n.$

In order to represent the temperature dependent material parameters $\tilde{k}$ and $\tilde{C}$ in a parametric form, we use the approximation properties of piecewise cubic interpolation methods, i.e. that for all $\epsilon>0$ there exists a number $n$ of partition points, such that 
\begin{align*}
    \sup\limits_{u \in U}|\tilde{k}(u)-k(u)|<\epsilon,\\
    \sup\limits_{u \in U}|\tilde{C}(u)-C(u)|<\epsilon.
\end{align*}
 
This means that, neglecting some small error, we can replace $\tilde{k}$ and $\tilde{C}$ by $k$ and $C$. Given the fixed partition $\pi_n$ and the interpolation method, the functions $k$ and $C$ can be represented by parameter vectors $(\underline{k},\underline{C})^T$.

With this we can define the {\em forward operator} explicitly by
\begin{align}
    \label{forwardOperator1}
    F:\ \mathbb{R}_+^{2n} &\to \mathscr{U},\\   \nonumber
     \underline{p}:=(\underline{k},\underline{C})^T&\mapsto u,
\end{align}
where $u\in \mathscr{U}$ is the solution to the modified initial boundary value problem (IBVP)
\begin{align}
\label{pde9}
C(u)u_t &= (k(u)u_z)_z, &(t,z)\in I \times \Omega \\\label{randUnten}
u&= u^b, &t\in I, z=0,\\ \label{randOben}
u&= u^t, &t\in I, z=L,\\ \label{pde12}
u&= u_0, &t=0, z\in \Omega,
\end{align}
to some material parameters represented by $(\underline{k},\underline{C})^T$ with known $u^b,u^t$ and $u_0.$

\subsection{The observation operator}
While a solution $u\in \mathcal{U}$ for the IBVP \eqref{pde9}-\eqref{pde12} is a function of time $t\in \bar{I}$ and position $z\in \bar{\Omega}$, our available data  $u^c$ consist of discrete temperature measurements at the core of the heavy plate, see Section \ref{sec:zwei}. That means the position $z=\frac{L}{2}$ is fixed and the temperature is recorded at times $t_j$ for $j=1,\dots,m,$ i.e. we have $u^c=(u_1^c,\dots,u_m^c)^T\in\mathbb{R}_+^m.$

To be able to compare the output of our forward operator with the data, we need an auxiliary operator, the so-called {\em observation operator}
\begin{align}
    \label{obsvervationOPERATOR}
    Q:\mathcal{U}&\to \mathbb{R}_+^m,\\
\nonumber    u&\mapsto \left(u\left(t_1,\frac{L}{2}\right),\dots,u\left(t_m,\frac{L}{2}\right)\right)^T.
\end{align}
This finally allows us to formulate the mathematical model of the inverse heat conduction problem.

\section{Parameter estimation of the material parameters $C$ and $k$}
\label{sec:vier}
\subsection{The formulation of the inverse heat conduction problem}

For given functions $u^b$, $u^t$, $u_0$ and some given data $u^c\in\mathbb{R}_+^m$, we want to deduce the heat conduction behaviour, i.e. find the optimal (interpolated) material parameters $k_{opt}$ and $C_{opt}$, represented by the function values $\underline{k}_{opt}$ and $\underline{C}_{opt}$ to the fixed partition $\pi_n:u_{min}=u_1<u_2<\cdots<u_n=u_{max}$, such that the observed temperature fits the data. More precisely, the goal is to determine
\begin{align} \label{koptCopt}
        \left(\underline{k}_{opt},\underline{C}_{opt}\right)^T=\argminA_{\underline{p}:=(\underline{k},\underline{C})^T\in \mathbb{R}^{2n}_+}\left\|QF(\underline{p})-u^c\right\|^2_2,
\end{align}
where $F$ is the forward operator \eqref{forwardOperator1}, $Q$ the observation operator \eqref{obsvervationOPERATOR} and $QF$ their composition. By $\|\cdot\|_2$ we denote the standard euclidean norm.

We can assume that given our data $u^c \in \mathbb{R}_+^m$ a solution \eqref{koptCopt} always exists.

Note, that $u^c \in \mathbb{R}^m_+$ contains small measurement errors in the sense that we have
\begin{align*}
    \sum\limits_{j=1}^m \left(u^{ex}\left(t_j,\frac{L}{2}\right)-u^c_j\right)^2\leq \delta
\end{align*}
for some $\delta>0$, where $u^{ex}\left(\cdot,\frac{L}{2}\right)$ denotes the time evaluation of the exact temperature $u^{ex}$ at the core position, which is of course not available to us in reality. Nevertheless, the noise level $\delta$ is considered small such that $u^c$ represents reliable data. Also, we plan to fully discretize the forward operator leading to a finite-dimensional range and thus to an inverse problem that is \textit{well-posed}, but probably \textit{ill-conditioned}, see \cite{louis}. This means, that the underlying inverse problem is not \textit{ill-posed} in the sense of Nashed (see \cite{schuster}) and at most mildly ill-conditioned. Hence, small errors in the data amplify the errors in the solution only slightly and a regularizing technique, e.g. by adding a penalty term in \eqref{koptCopt}, is not necessary. We recommend \cite{schuster} for readers who are interested in regularization methods in general settings.

In the following subsection, we want to make an important remark about the ambiguity of solutions.

\subsection{Ambiguous solution of the inverse problem}\label{subsection:vierzwei}

We would like to mention that the solution pair $ \left(\underline{k}_{opt},\underline{C}_{opt}\right)^T\in \mathbb{R}^{2n}_+$ can't be the unique minimizer of the least-squares functional
\begin{equation}
    \label{functional}
    J(\underline{p}):=\left\|QF(\underline{p})-u^c\right\|^2_2, 
\end{equation}
because the forward operator $F$ is not injective, i.e., for any real number $\alpha>0$ we have

\begin{equation}
    \label{FnichtInjektiv}
    F(\underline{p})=F(\alpha \underline{p})\ \ \  \forall\ \underline{p} \in \mathbb{R}_+^{2n}.
\end{equation}

This follows from the fact that multiplying equation \eqref{pde9} by some scalar $\alpha$ doesn't influence the differential equation. Thus, it is impossible to eliminate the ambiguity of the inverse problem solution. Rather, one has to accept the solution, knowing that only the quotient 
\begin{equation}
    \lambda_{opt}(u):=\frac{k_{opt}(u)}{C_{opt}(u)}.
\end{equation}
is uniquely determined.

This way, we can deal with the non-injectivity \eqref{FnichtInjektiv} of the forward operator, because $\underline{p}=(\underline{k},\underline{C})^T$ and $\alpha \underline{p}=(\alpha\underline{k},\alpha\underline{C})^T$ are from the same class sharing the same quotient of the corresponding interpolated functions, i.e.
\begin{equation}
    \frac{k(u)}{C(u)}=\lambda(u)=\frac{\alpha k(u)}{\alpha C(u)}.
\end{equation}

In thermodynamics, $\lambda$ is the so-called {\em thermal diffusivity} of the material.

Here, the thermal conduction behaviour of the material depends mainly on the thermal diffusivity $\lambda$. A perfectly decoupled characterization of $k$ and $C$ is not possible.
Interestingly, this coincides with the insights for the solution of \textit{linear} inverse heat conduction problems where only the scalar thermal diffusivity is discussed. 

\underline{Remark:} 
In order to be able to identify a heat flux $q(t)$ with Fourier's law of heat conduction in the form of
\begin{align*}
    \pm k(u)u_z = q
\end{align*}
on the surfaces of the heavy plate, we need to know the functional form of $k(u)$.
However, the identification of some solution pair $(\alpha \underline{k},\alpha \underline{C})^T$ is good enough, if one accepts that the determination of $q$ is only disturbed by a scalar and time-independent factor $\alpha$, i.e.
\begin{align*}
    \pm \alpha k(u)u_z &= \tilde{q},\\
    q &= \frac{\tilde{q}}{\alpha}.
\end{align*} 
This is the main reason why we propose a parameter estimation problem to simultaneously determine both, the functional forms of $k$ and $C$. 

In Section \ref{sec:fuenf} we discuss the numerical implementation of solving the parameter estimation problem \eqref{koptCopt}.

\section{Implementation approach and numerical results} \label{sec:fuenf}

While the previous sections were more introductory and theoretical, in this section we want to give a practical guide on how to numerically solve the parameter estimation problem, i.e. the minimization 
\begin{equation}
    \label{functionEstimationMinmizationProblem}
    \min\limits_{\underline{p}\in \mathbb{R}_+^{2n}}J(\underline{p})=\min\limits_{\underline{p}\in \mathbb{R}_+^{2n}}\left\|QF(\underline{p})-u^c\right\|^2_2.
\end{equation}

Assuming that we recorded $M\geq 1$ experiments, we can even generalize the objective functional to
\begin{equation}
    \label{MExperiments}
    J_M(\underline{p})=\sum\limits_{i=1}^M\left\|QF_i(\underline{p})-u^{c,i}\right\|^2_2,
\end{equation}
where $F_i$ maps the parameter vector $\underline{p}=(\underline{k},\underline{C})^T$ to the solution of the initial boundary value problem
\begin{align}
\label{pdeM_Anfang}
C(u)u_t &= (k(u)u_z)_z, &(t,z)\in I \times \Omega \\
u&= u^{b,i}, &t\in I, z=0,\\
u&= u^{t,i}, &t\in I, z=L,\\ \label{pdeM_Ende}
u&= u_0^i, &t=0, z\in \Omega,
\end{align}
for $i=1,\dots,M.$ Note, that the minimization \eqref{functionEstimationMinmizationProblem} represents the special case $M=1.$ Recording and using several experiments stabilize the minimization process and yield better results in case of noisy core measurements $u^{c,i}$. \\ \ \\

In this section, we address the following topics:

\begin{itemize}
    \item[(a)] Representation of temperature dependent material parameters $k(u)$, $C(u) \in \mathcal{C}^1(U)$ by some parameter vector $\underline{p}=(\underline{k},\underline{C})^T \in \mathbb{R}_+^{2n}$. 
    \item[(b)] Implementation of the forward operator $F_i$, mapping a vector $\underline{p}$ to some temperature matrix, i.e. some time and space discretization of $u$.
    \item[(c)] Application of the observation operator to the solution of the forward problem.
    \item[(d)] Solving \eqref{MExperiments} for \textit{simulated} data $u^{c,i}_{sim}$ ($i=1,\dots,M$) and comparison of the simulated material parameters $k_{sim}$ and $C_{sim}$ to the optimized functions $k_{opt}$ and $C_{opt}$. 
\end{itemize}

\subsection*{(a):}

The following paragraph only deals with the representation of the thermal conductivity $k(u)$ by some parameter vector $\underline{k} \in \mathbb{R}^n_+$, but the procedure works in the same way for the volumetric heat capacity $C(u).$

First of all, we want to fix the number $n>0$ of partition points $(u_1,\dots,u_n)^T \in \mathbb{R}^n$ of the interval $U=[u_{min},u_{max}]$. By also setting the parameter vector $\underline{k}=(k_1,\dots,k_n)^T\in \mathbb{R}^n_+$ of the corresponding function values and choosing an interpolation method, we can create a function $k(u) \in \mathcal{C}^1(U)$ by interpolation, such that 
\begin{align*}
    k(u_i)=k_i
\end{align*}
for $i=1,\dots,n$.

In our case, we want to use the Piecewise Cubic Hermite Interpolating Polynomials (PCHIP) introduced in \cite{fritschCarlson}, due to the favorable monotonic behaviour of the interpolant with respect to the function values. Of course, one could also use an alternative interpolation method.

A big advantage in using an interpolation method rather than some representation as a linear combination of functions from a physically meaningful dictionary is that no a-priori information about the functional form of the material parameters are needed. Also, adjusting one parameter $k_i$ of course leads to some changed interpolant $k(u)$, but only the function values in the local region of the corresponding $u_i$ are affected (see Figure \ref{fig:pchip}), which is very useful in the minimization process. 
\begin{align}
\label{dreiundzwanzig}
    \min\limits_{\underline{p}\in \mathbb{R}_+^{2n}}J_M(\underline{p}).
\end{align}

\begin{center}
\def\svgwidth{350pt}
\begingroup%
  \makeatletter%
  \providecommand\color[2][]{%
    \errmessage{(Inkscape) Color is used for the text in Inkscape, but the package 'color.sty' is not loaded}%
    \renewcommand\color[2][]{}%
  }%
  \providecommand\transparent[1]{%
    \errmessage{(Inkscape) Transparency is used (non-zero) for the text in Inkscape, but the package 'transparent.sty' is not loaded}%
    \renewcommand\transparent[1]{}%
  }%
  \providecommand\rotatebox[2]{#2}%
  \newcommand*\fsize{\dimexpr\f@size pt\relax}%
  \newcommand*\lineheight[1]{\fontsize{\fsize}{#1\fsize}\selectfont}%
  \ifx\svgwidth\undefined%
    \setlength{\unitlength}{1440bp}%
    \ifx\svgscale\undefined%
      \relax%
    \else%
      \setlength{\unitlength}{\unitlength * \real{\svgscale}}%
    \fi%
  \else%
    \setlength{\unitlength}{\svgwidth}%
  \fi%
  \global\let\svgwidth\undefined%
  \global\let\svgscale\undefined%
  \makeatother%
  \begin{picture}(1,0.4234375)%
    \lineheight{1}%
    \setlength\tabcolsep{0pt}%
    \put(0,0){\includegraphics[width=\unitlength,page=1]{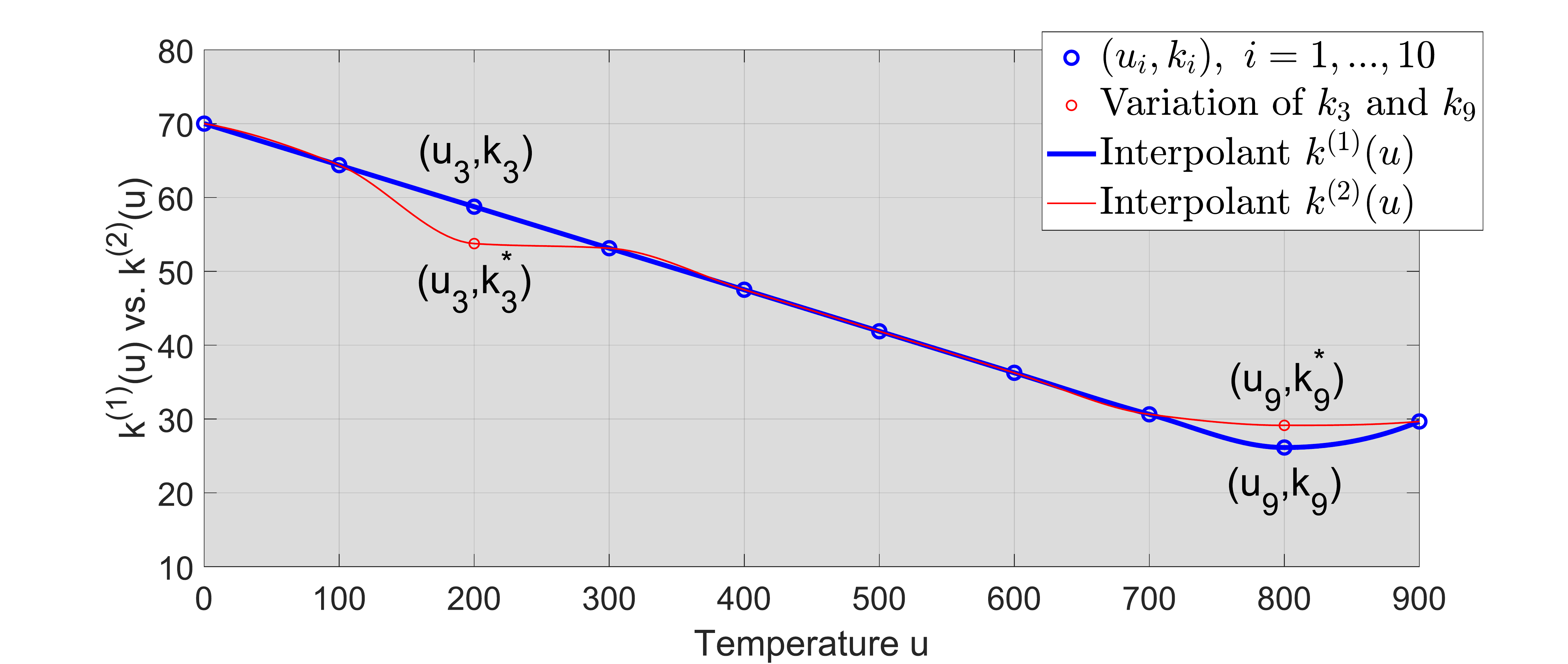}}%
  \end{picture}%
\endgroup%

\captionof{figure}{Examples of PCHIP interpolants for $U=[0,900]$, $u_i=\frac{900\cdot (i-1)}{(n-1)}$, $i=1,\dots,n=10$, with varying $k_3$ and $k_9$}
\label{fig:pchip}
\end{center}

Increasing the number $n$ of partition points allows us to represent functions that are more complex. Thus, for the estimation of the unknown material parameters $k(u)$ and $C(u)$, we want to fix a sufficiently large number $n$ of partition points to determine enough function values $(\underline{k},\underline{C})^T$ which fit the model and the observed data in anoptimal way.

\subsection*{(b):}

In this subsection we present the implementation of the forward operator $F$, which maps the parameter vector of function values $(\underline{k},\underline{C})^T$ to a solution $u\in \mathcal{U}$ of the initial boundary value problem \eqref{pde9}-\eqref{pde12}. The forward operators $F_i$ are implemented analogously.

To get a solution of the IBVP numerically, we discretize the temperature  $u \in \mathcal{U}$ with respect to time $t\in I=[0,T]$ and the space variable $z\in\Omega=[0,L]$. 
For that reason, we define equidistant partitions
\begin{align}\label{partitionT}
    0=t_1<t_2<\dots<t_m=T,\\
    0=z_1<z_2<\dots<z_l=L \label{partitonL}
\end{align}
for $I$ and $\Omega$ with increments $\Delta t=\frac{T}{m-1}$ and $\Delta z=\frac{L}{l-1}$, respectively. Here, $m,l>0$ are sufficiently large integers. 


\begin{center}
\def\svgwidth{350pt}
\begingroup%
  \makeatletter%
  \providecommand\color[2][]{%
    \errmessage{(Inkscape) Color is used for the text in Inkscape, but the package 'color.sty' is not loaded}%
    \renewcommand\color[2][]{}%
  }%
  \providecommand\transparent[1]{%
    \errmessage{(Inkscape) Transparency is used (non-zero) for the text in Inkscape, but the package 'transparent.sty' is not loaded}%
    \renewcommand\transparent[1]{}%
  }%
  \providecommand\rotatebox[2]{#2}%
  \newcommand*\fsize{\dimexpr\f@size pt\relax}%
  \newcommand*\lineheight[1]{\fontsize{\fsize}{#1\fsize}\selectfont}%
  \ifx\svgwidth\undefined%
    \setlength{\unitlength}{501.73664312bp}%
    \ifx\svgscale\undefined%
      \relax%
    \else%
      \setlength{\unitlength}{\unitlength * \real{\svgscale}}%
    \fi%
  \else%
    \setlength{\unitlength}{\svgwidth}%
  \fi%
  \global\let\svgwidth\undefined%
  \global\let\svgscale\undefined%
  \makeatother%
  \begin{picture}(1,0.62842282)%
    \lineheight{1}%
    \setlength\tabcolsep{0pt}%
    \put(0,0){\includegraphics[width=\unitlength,page=1]{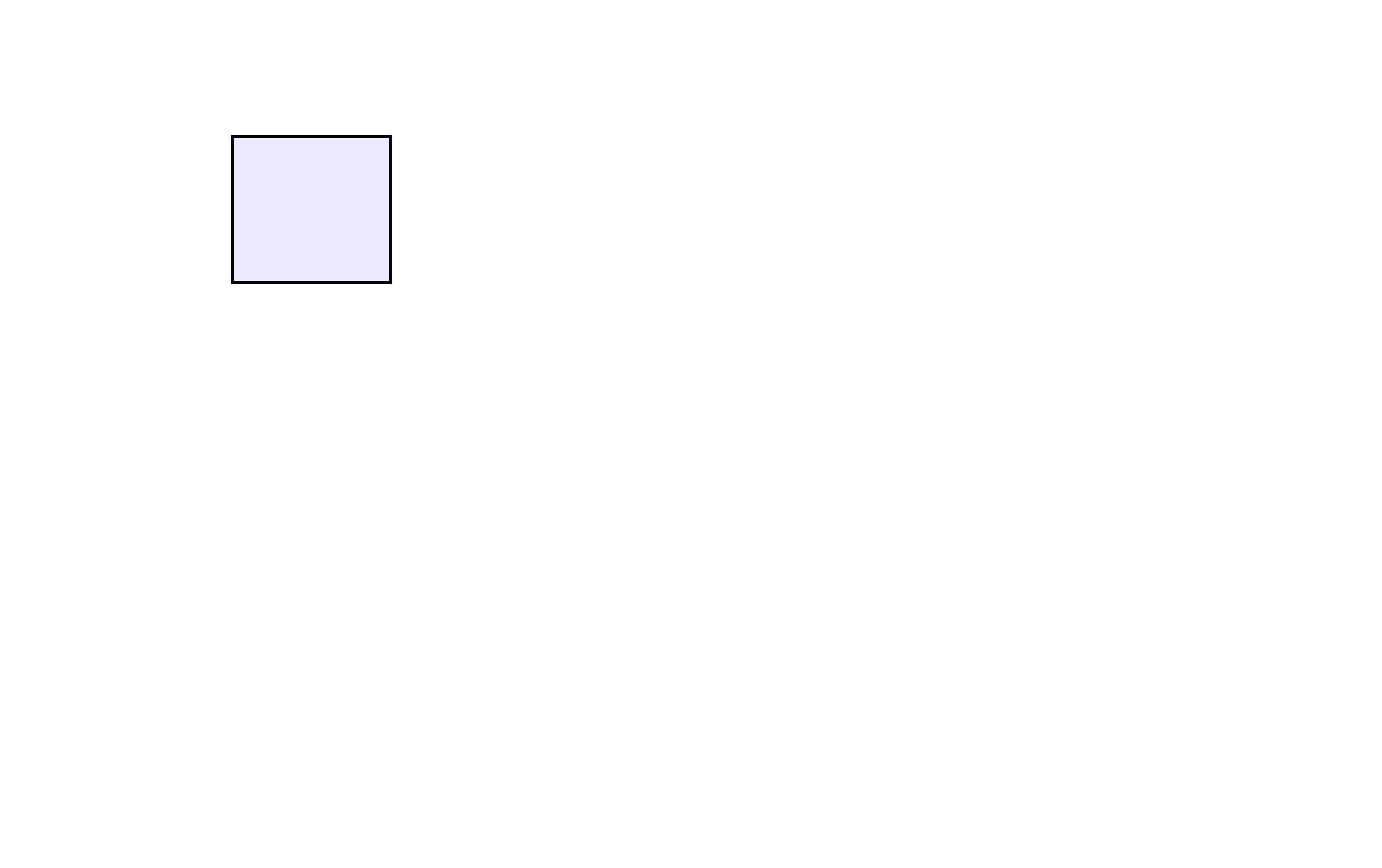}}%
    \put(0.20103116,0.46916577){\color[rgb]{0,0,0}\makebox(0,0)[lt]{\lineheight{1.25}\smash{\begin{tabular}[t]{l}$u_{11}$\end{tabular}}}}%
    \put(0,0){\includegraphics[width=\unitlength,page=2]{tempMatrix.pdf}}%
    \put(0.31186263,0.46944448){\color[rgb]{0,0,0}\makebox(0,0)[lt]{\lineheight{1.25}\smash{\begin{tabular}[t]{l}$u_{12}$\end{tabular}}}}%
    \put(0.42816174,0.46944448){\color[rgb]{0,0,0}\makebox(0,0)[lt]{\lineheight{1.25}\smash{\begin{tabular}[t]{l}$u_{13}$\end{tabular}}}}%
    \put(0.63882357,0.46938379){\color[rgb]{0,0,0}\makebox(0,0)[lt]{\lineheight{1.25}\smash{\begin{tabular}[t]{l}$u_{1,l-2}$\end{tabular}}}}%
    \put(0.75541837,0.46938379){\color[rgb]{0,0,0}\makebox(0,0)[lt]{\lineheight{1.25}\smash{\begin{tabular}[t]{l}$u_{1,l-1}$\end{tabular}}}}%
    \put(0.88695971,0.46938379){\color[rgb]{0,0,0}\makebox(0,0)[lt]{\lineheight{1.25}\smash{\begin{tabular}[t]{l}$u_{1l}$\end{tabular}}}}%
    \put(0.20036564,0.36463811){\color[rgb]{0,0,0}\makebox(0,0)[lt]{\lineheight{1.25}\smash{\begin{tabular}[t]{l}$u_{21}$\end{tabular}}}}%
    \put(0.31119728,0.36491682){\color[rgb]{0,0,0}\makebox(0,0)[lt]{\lineheight{1.25}\smash{\begin{tabular}[t]{l}$u_{22}$\end{tabular}}}}%
    \put(0.75475776,0.36485613){\color[rgb]{0,0,0}\makebox(0,0)[lt]{\lineheight{1.25}\smash{\begin{tabular}[t]{l}$u_{2,l-1}$\end{tabular}}}}%
    \put(0.88629911,0.36485613){\color[rgb]{0,0,0}\makebox(0,0)[lt]{\lineheight{1.25}\smash{\begin{tabular}[t]{l}$u_{2l}$\end{tabular}}}}%
    \put(0.18242793,0.15237518){\color[rgb]{0,0,0}\makebox(0,0)[lt]{\lineheight{1.25}\smash{\begin{tabular}[t]{l}$u_{m-1,1}$\end{tabular}}}}%
    \put(0.86836165,0.1525932){\color[rgb]{0,0,0}\makebox(0,0)[lt]{\lineheight{1.25}\smash{\begin{tabular}[t]{l}$u_{m-1,l}$\end{tabular}}}}%
    \put(0.20036319,0.04773864){\color[rgb]{0,0,0}\makebox(0,0)[lt]{\lineheight{1.25}\smash{\begin{tabular}[t]{l}$u_{m1}$\end{tabular}}}}%
    \put(0.311195,0.04801744){\color[rgb]{0,0,0}\makebox(0,0)[lt]{\lineheight{1.25}\smash{\begin{tabular}[t]{l}$u_{m2}$\end{tabular}}}}%
    \put(0.75476026,0.0479567){\color[rgb]{0,0,0}\makebox(0,0)[lt]{\lineheight{1.25}\smash{\begin{tabular}[t]{l}$u_{m,l-1}$\end{tabular}}}}%
    \put(0.88630161,0.0479567){\color[rgb]{0,0,0}\makebox(0,0)[lt]{\lineheight{1.25}\smash{\begin{tabular}[t]{l}$u_{ml}$\end{tabular}}}}%
    \put(0,0){\includegraphics[width=\unitlength,page=3]{tempMatrix.pdf}}%
    \put(-0.00150843,0.4679345){\color[rgb]{0,0,0}\makebox(0,0)[lt]{\lineheight{1.25}\smash{\begin{tabular}[t]{l}$t=t_1$\end{tabular}}}}%
    \put(-0.00123107,0.36598486){\color[rgb]{0,0,0}\makebox(0,0)[lt]{\lineheight{1.25}\smash{\begin{tabular}[t]{l}$t=t_2$\end{tabular}}}}%
    \put(-0.00117035,0.04820528){\color[rgb]{0,0,0}\makebox(0,0)[lt]{\lineheight{1.25}\smash{\begin{tabular}[t]{l}$t=t_m$\end{tabular}}}}%
    \put(0.1682429,0.61280676){\color[rgb]{0,0,0}\makebox(0,0)[lt]{\lineheight{1.25}\smash{\begin{tabular}[t]{l}$z=z_1$\end{tabular}}}}%
    \put(0.28804322,0.61308546){\color[rgb]{0,0,0}\makebox(0,0)[lt]{\lineheight{1.25}\smash{\begin{tabular}[t]{l}$z=z_2$\end{tabular}}}}%
    \put(0.40733247,0.61308546){\color[rgb]{0,0,0}\makebox(0,0)[lt]{\lineheight{1.25}\smash{\begin{tabular}[t]{l}$z=z_3$\end{tabular}}}}%
    \put(0.86314552,0.61302478){\color[rgb]{0,0,0}\makebox(0,0)[lt]{\lineheight{1.25}\smash{\begin{tabular}[t]{l}$z=z_l$\end{tabular}}}}%
    \put(0,0){\includegraphics[width=\unitlength,page=4]{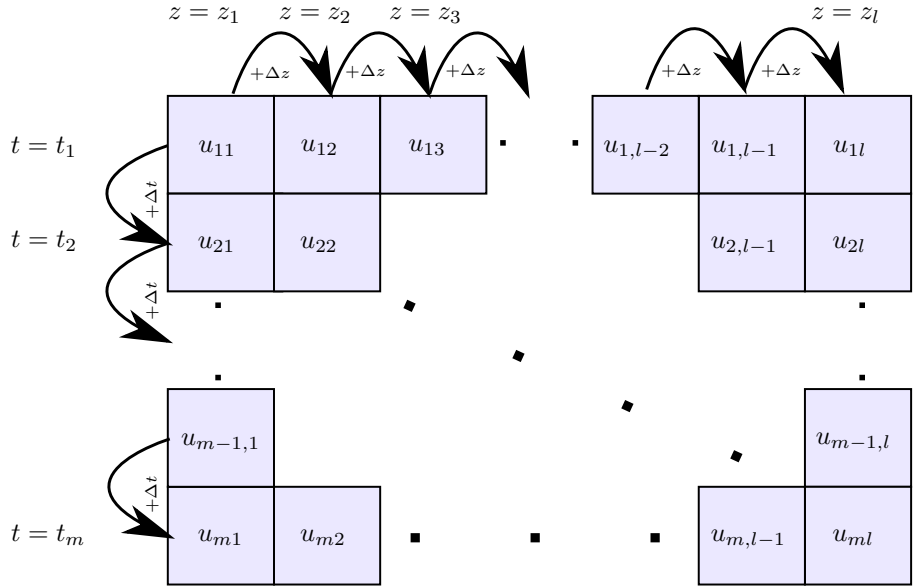}}%
    \put(0.25672035,0.5493228){\color[rgb]{0,0,0}\makebox(0,0)[lt]{\lineheight{1.25}\smash{\begin{tabular}[t]{l}\tiny{$+\Delta z$}\end{tabular}}}}%
    \put(0.35980621,0.55083279){\color[rgb]{0,0,0}\makebox(0,0)[lt]{\lineheight{1.25}\smash{\begin{tabular}[t]{l}\tiny{$+\Delta z$}\end{tabular}}}}%
    \put(0.46810015,0.55209998){\color[rgb]{0,0,0}\makebox(0,0)[lt]{\lineheight{1.25}\smash{\begin{tabular}[t]{l}\tiny{$+\Delta z$}\end{tabular}}}}%
    \put(0.70129031,0.55209998){\color[rgb]{0,0,0}\makebox(0,0)[lt]{\lineheight{1.25}\smash{\begin{tabular}[t]{l}\tiny{$+\Delta z$}\end{tabular}}}}%
    \put(0.80891565,0.55209998){\color[rgb]{0,0,0}\makebox(0,0)[lt]{\lineheight{1.25}\smash{\begin{tabular}[t]{l}\tiny{$+\Delta z$}\end{tabular}}}}%
    \put(0.1557859,0.39658687){\color[rgb]{0,0,0}\rotatebox{90}{\makebox(0,0)[lt]{\lineheight{1.25}\smash{\begin{tabular}[t]{l}\tiny{$+\Delta t$}\end{tabular}}}}}%
    \put(0.15581623,0.28724468){\color[rgb]{0,0,0}\rotatebox{90}{\makebox(0,0)[lt]{\lineheight{1.25}\smash{\begin{tabular}[t]{l}\tiny{$+\Delta t$}\end{tabular}}}}}%
    \put(0.15581623,0.0779713){\color[rgb]{0,0,0}\rotatebox{90}{\makebox(0,0)[lt]{\lineheight{1.25}\smash{\begin{tabular}[t]{l}\tiny{$+\Delta t$}\end{tabular}}}}}%
  \end{picture}%
\endgroup%

\captionof{figure}{Illustration of the discretization of $u$}
\label{fig:tempMatrix}
\end{center}

Adopting this notation, we can interpret the output of $F$ as a temperature matrix $u \in \mathbb{R}_+^{m\times l}$ with components
\begin{align}
    u_{ij}:=u(t_i,z_j),\ \ \ \  &i=1,\dots,m; \\
    &j=1,\dots,l,
\end{align}
see Figure \ref{fig:tempMatrix}.

However, these components are yet to be determined from the IBVP. From the boundary conditions \eqref{randUnten}-\eqref{randOben} and the intial temperature distribution \eqref{pde12} we get
\begin{align*}
    u_{i1}=u^b(t_i),\ \ \  & i=1,\dots m \\
    u_{il}=u^t(t_i),\ \ \ & i=1,\dots m \\
    u_{1j}=u_0(z_j),\ \ & j=1,\dots,l
\end{align*}
respectively. Finally, we determine
\begin{align*}
    u_{ij} \text{\ \  for \ \ } i=2,\dots,m \text{\ \  and \ \ } j=2,\dots l-1
\end{align*}
from \eqref{pde9} by using a finite difference method which is the following marching scheme 
\begin{align}
    \label{marchingScheme}
    u_{ij}=u_{i-1,j}+\frac{\Delta t}{(\Delta z) ^2 C(u_{i-1,j})}\left(k_{ij}^e\cdot(u_{i-1,j+1}-u_{i-1,j})-k_{ij}^w\cdot(u_{i-1,j}-u_{i-1,j-1})\right),
\end{align}
cf. Figure \ref{fig:tempMatrix}, with the harmonic means of thermal conductivities $$k_{ij}^e:=\frac{2\cdot k(u_{i-1,j+1})\cdot k(u_{i-1,j})}{k(u_{i-1,j+1})+k(u_{i-1,j})},$$ \ \\
$$k_{ij}^w:=\frac{2\cdot k(u_{i-1,j-1})\cdot k(u_{i-1,j})}{k(u_{i-1,j-1})+k(u_{i-1,j})}.$$
Note, that we choose an explicit scheme, i.e. the new time step component $u_{ij}$ is an explicit function of the old time step components $u_{i-1,j-1},\ u_{i-1,j}$ and $u_{i-1,j+1}$. This way, by ensuring that the time and space discretizations are chosen carefully due to instability issues, the implementation of \eqref{marchingScheme} is straightforward. 
We achieve overall good results by guaranteeing the stability condition 
\begin{align}
\label{stabilitaetsvoraussetzung}
\Delta t \leq \frac{(\Delta z)^2}{2\max\limits_{\tilde{u}\in U}\lambda(\tilde{u})}.
\end{align}

\subsection*{(c):}

Now that we interpret the output of $F$ as a temperature matrix $u \in \mathbb{R}_+^{m\times l}$, we modify the observation operator $Q$ to be
\begin{align}
    \label{obsvervationOPERATORmatrix}
    Q:\mathbb{R}_+^{m\times l}&\to \mathbb{R}_+^m,\\
\nonumber    u&\mapsto \left(u\left(t_1,\frac{L}{2}\right),\dots,u\left(t_m,\frac{L}{2}\right)\right)^T.
\end{align}
Here, we choose the partition of the time interval \eqref{partitionT} to fit the data measurement setting, i.e. $u^c \in \mathbb{R}^m_+$ and $u^c_i$ is the core temperature measured at times $t_i$ for $i=1,\dots,m.$ Also, we require the partition of the space interval \eqref{partitonL} to contain the depth $z_k$ for some $1<k<l$ with $z_k=\frac{L}{2},$ i.e. the space discretization yields a depth corresponding to the core of the thickness dimension of the heavy plate. Nevertheless, we like to mention that it is possible to choose the partitions \eqref{partitionT}-\eqref{partitonL} in every possible way as long as the marching scheme \eqref{marchingScheme} is numerically stable, cf. \eqref{stabilitaetsvoraussetzung}. In this case however, it will be probably necessary to interpolate the temperature matrix $u$ such that we can extract a temperature vector $$\left(u\left(t_1,\frac{L}{2}\right),\dots,u\left(t_m,\frac{L}{2}\right)\right)^T.$$
\begin{center}
\def\svgwidth{400pt}
    \footnotesize{
\begingroup%
  \makeatletter%
  \providecommand\color[2][]{%
    \errmessage{(Inkscape) Color is used for the text in Inkscape, but the package 'color.sty' is not loaded}%
    \renewcommand\color[2][]{}%
  }%
  \providecommand\transparent[1]{%
    \errmessage{(Inkscape) Transparency is used (non-zero) for the text in Inkscape, but the package 'transparent.sty' is not loaded}%
    \renewcommand\transparent[1]{}%
  }%
  \providecommand\rotatebox[2]{#2}%
  \newcommand*\fsize{\dimexpr\f@size pt\relax}%
  \newcommand*\lineheight[1]{\fontsize{\fsize}{#1\fsize}\selectfont}%
  \ifx\svgwidth\undefined%
    \setlength{\unitlength}{1397.2499824bp}%
    \ifx\svgscale\undefined%
      \relax%
    \else%
      \setlength{\unitlength}{\unitlength * \real{\svgscale}}%
    \fi%
  \else%
    \setlength{\unitlength}{\svgwidth}%
  \fi%
  \global\let\svgwidth\undefined%
  \global\let\svgscale\undefined%
  \makeatother%
  \begin{picture}(1,0.50295223)%
    \lineheight{1}%
    \setlength\tabcolsep{0pt}%
    \put(0,0){\includegraphics[width=\unitlength,page=1]{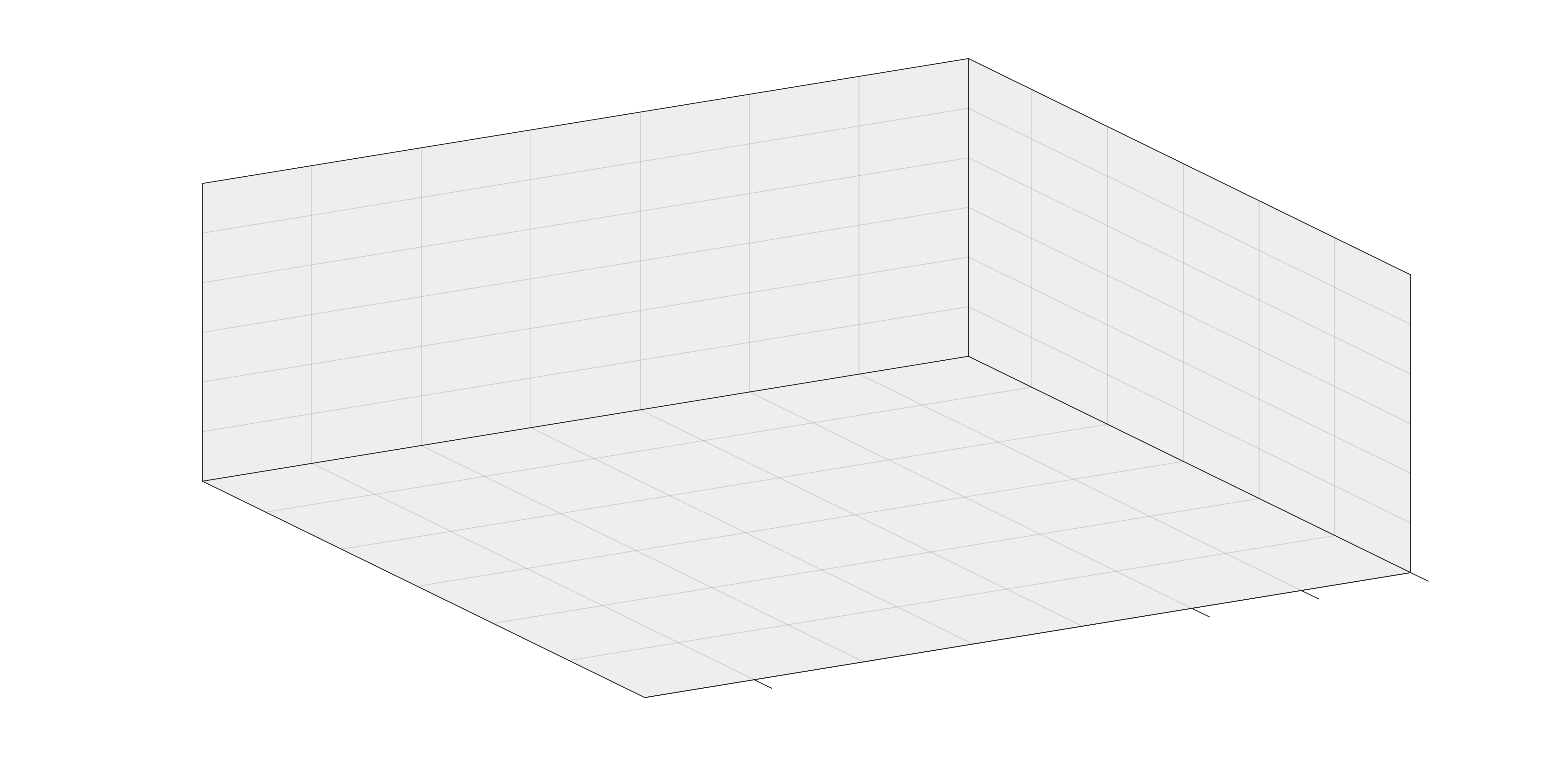}}%
    \put(0.84905899,0.10735674){\makebox(0,0)[lt]{\lineheight{1.25}\smash{\begin{tabular}[t]{l}60\end{tabular}}}}%
    \put(0,0){\includegraphics[width=\unitlength,page=2]{observationOperator.pdf}}%
    \put(0.77886044,0.09590456){\makebox(0,0)[lt]{\lineheight{1.25}\smash{\begin{tabular}[t]{l}50\end{tabular}}}}%
    \put(0,0){\includegraphics[width=\unitlength,page=3]{observationOperator.pdf}}%
    \put(0.70866189,0.08445233){\makebox(0,0)[lt]{\lineheight{1.25}\smash{\begin{tabular}[t]{l}40\end{tabular}}}}%
    \put(0,0){\includegraphics[width=\unitlength,page=4]{observationOperator.pdf}}%
    \put(0.57648658,0.0333591){\makebox(0,0)[lt]{\lineheight{1.25}\smash{\begin{tabular}[t]{l}Time t in seconds\end{tabular}}}}%
    \put(0,0){\includegraphics[width=\unitlength,page=5]{observationOperator.pdf}}%
    \put(0.63846334,0.07300011){\makebox(0,0)[lt]{\lineheight{1.25}\smash{\begin{tabular}[t]{l}30\end{tabular}}}}%
    \put(0,0){\includegraphics[width=\unitlength,page=6]{observationOperator.pdf}}%
    \put(0.56826473,0.06154793){\makebox(0,0)[lt]{\lineheight{1.25}\smash{\begin{tabular}[t]{l}20\end{tabular}}}}%
    \put(0,0){\includegraphics[width=\unitlength,page=7]{observationOperator.pdf}}%
    \put(0.49806618,0.05009565){\makebox(0,0)[lt]{\lineheight{1.25}\smash{\begin{tabular}[t]{l}10\end{tabular}}}}%
    \put(0,0){\includegraphics[width=\unitlength,page=8]{observationOperator.pdf}}%
    \put(0.42786763,0.03864348){\makebox(0,0)[lt]{\lineheight{1.25}\smash{\begin{tabular}[t]{l}0\end{tabular}}}}%
    \put(0,0){\includegraphics[width=\unitlength,page=9]{observationOperator.pdf}}%
    \put(0.34478674,0.04199452){\makebox(0,0)[lt]{\lineheight{1.25}\smash{\begin{tabular}[t]{l}$z=0$\end{tabular}}}}%
    \put(0,0){\includegraphics[width=\unitlength,page=10]{observationOperator.pdf}}%
    \put(0.1798175,0.11433172){\makebox(0,0)[lt]{\lineheight{1.25}\smash{\begin{tabular}[t]{l}$z=\frac{L}{2}$\end{tabular}}}}%
    \put(0,0){\includegraphics[width=\unitlength,page=11]{observationOperator.pdf}}%
    \put(0.1354745,0.08248846){\makebox(0,0)[lt]{\lineheight{1.25}\smash{\begin{tabular}[t]{l}Depth z\end{tabular}}}}%
    \put(0,0){\includegraphics[width=\unitlength,page=12]{observationOperator.pdf}}%
    \put(0.08503661,0.39192641){\makebox(0,0)[lt]{\lineheight{1.25}\smash{\begin{tabular}[t]{l}800\end{tabular}}}}%
    \put(0.08503661,0.36008057){\makebox(0,0)[lt]{\lineheight{1.25}\smash{\begin{tabular}[t]{l}700\end{tabular}}}}%
    \put(0.08503661,0.32823473){\makebox(0,0)[lt]{\lineheight{1.25}\smash{\begin{tabular}[t]{l}600\end{tabular}}}}%
    \put(0.08503661,0.29638889){\makebox(0,0)[lt]{\lineheight{1.25}\smash{\begin{tabular}[t]{l}500\end{tabular}}}}%
    \put(0.08503661,0.26454305){\makebox(0,0)[lt]{\lineheight{1.25}\smash{\begin{tabular}[t]{l}400\end{tabular}}}}%
    \put(0.08503661,0.23269721){\makebox(0,0)[lt]{\lineheight{1.25}\smash{\begin{tabular}[t]{l}300\end{tabular}}}}%
    \put(0.08503661,0.20085137){\makebox(0,0)[lt]{\lineheight{1.25}\smash{\begin{tabular}[t]{l}200\end{tabular}}}}%
    \put(0.06967096,0.17547756){\makebox(0,0)[lt]{\lineheight{1.25}\smash{\begin{tabular}[t]{l}$z=L$\end{tabular}}}}%
    \put(0.07214761,0.22553167){\rotatebox{90}{\makebox(0,0)[lt]{\lineheight{1.25}\smash{\begin{tabular}[t]{l}Temperature in °C\end{tabular}}}}}%
    \put(0,0){\includegraphics[width=\unitlength,page=13]{observationOperator.pdf}}%
    \put(0.58346753,0.424584){\makebox(0,0)[lt]{\lineheight{1.25}\smash{\begin{tabular}[t]{l}Observation of core temp.\end{tabular}}}}%
    \put(0,0){\includegraphics[width=\unitlength,page=14]{observationOperator.pdf}}%
    \put(0.96136694,0.07103038){\makebox(0,0)[lt]{\lineheight{1.25}\smash{\begin{tabular}[t]{l}250\end{tabular}}}}%
    \put(0.96136694,0.10712507){\makebox(0,0)[lt]{\lineheight{1.25}\smash{\begin{tabular}[t]{l}300\end{tabular}}}}%
    \put(0.96136694,0.14321975){\makebox(0,0)[lt]{\lineheight{1.25}\smash{\begin{tabular}[t]{l}350\end{tabular}}}}%
    \put(0.96136694,0.17931449){\makebox(0,0)[lt]{\lineheight{1.25}\smash{\begin{tabular}[t]{l}400\end{tabular}}}}%
    \put(0.96136694,0.21540918){\makebox(0,0)[lt]{\lineheight{1.25}\smash{\begin{tabular}[t]{l}450\end{tabular}}}}%
    \put(0.96136694,0.25150386){\makebox(0,0)[lt]{\lineheight{1.25}\smash{\begin{tabular}[t]{l}500\end{tabular}}}}%
    \put(0.96136694,0.28759855){\makebox(0,0)[lt]{\lineheight{1.25}\smash{\begin{tabular}[t]{l}550\end{tabular}}}}%
    \put(0.96136694,0.32369329){\makebox(0,0)[lt]{\lineheight{1.25}\smash{\begin{tabular}[t]{l}600\end{tabular}}}}%
    \put(0.96136694,0.35978798){\makebox(0,0)[lt]{\lineheight{1.25}\smash{\begin{tabular}[t]{l}650\end{tabular}}}}%
    \put(0.96136694,0.39588266){\makebox(0,0)[lt]{\lineheight{1.25}\smash{\begin{tabular}[t]{l}700\end{tabular}}}}%
    \put(0.96136694,0.43197735){\makebox(0,0)[lt]{\lineheight{1.25}\smash{\begin{tabular}[t]{l}750\end{tabular}}}}%
    \put(0,0){\includegraphics[width=\unitlength,page=15]{observationOperator.pdf}}%
  \end{picture}%
\endgroup%
}

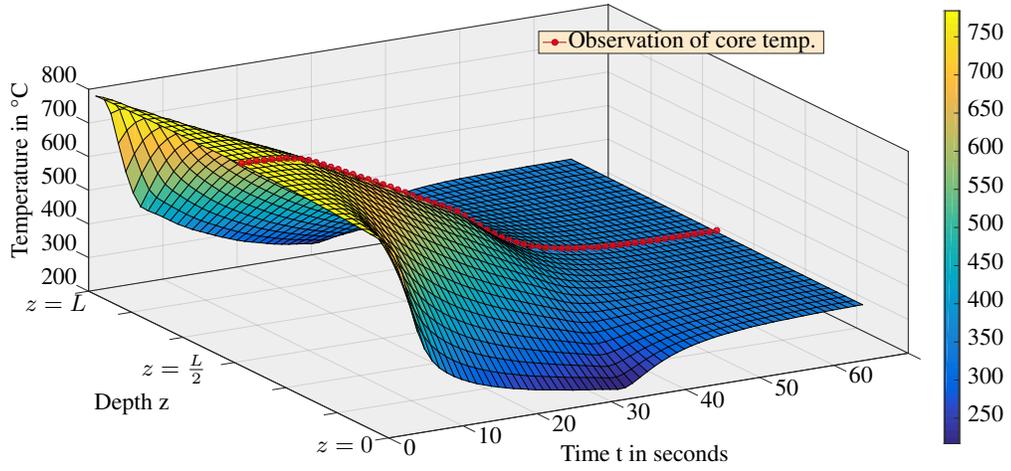
\captionof{figure}{Observing the core temperature at times $t_i$ from an interpolated temperature matrix $u$}
\label{fig:observationOperator}
\end{center}

In this way the observation operator can be interpreted as that it 'looks' at the temperature matrix and 'picks' the temperature corresponding to the core temperature at given time instances $t_i$ ($i=1,\dots,m),$ see also Figure \ref{fig:observationOperator}.

\subsection*{(d):}

In this part, we demonstrate the reliability of the implementation of \eqref{dreiundzwanzig} by means of simulations with synthetic data. Let $k_{sim}$ and $C_{sim}$ be some known (physically plausible) temperature dependent material parameters. In particular, for temperatures $u \in U=[0,900]$ we define

\begin{align}
    k_{sim}(u) &= 60-\frac{u}{30}, \\
    C_{sim}(u) &= 7650*\left(475+0.0265u+0.000855u^2-\frac{0.000855u^2-0.1735u+140}{1+e^{-0.1(u-700)}}\right),
\end{align}
see Figure \ref{fig:ksimCsim}.

\begin{center}
\def\svgwidth{450pt}
\begingroup%
  \makeatletter%
  \providecommand\color[2][]{%
    \errmessage{(Inkscape) Color is used for the text in Inkscape, but the package 'color.sty' is not loaded}%
    \renewcommand\color[2][]{}%
  }%
  \providecommand\transparent[1]{%
    \errmessage{(Inkscape) Transparency is used (non-zero) for the text in Inkscape, but the package 'transparent.sty' is not loaded}%
    \renewcommand\transparent[1]{}%
  }%
  \providecommand\rotatebox[2]{#2}%
  \newcommand*\fsize{\dimexpr\f@size pt\relax}%
  \newcommand*\lineheight[1]{\fontsize{\fsize}{#1\fsize}\selectfont}%
  \ifx\svgwidth\undefined%
    \setlength{\unitlength}{1222.6206665bp}%
    \ifx\svgscale\undefined%
      \relax%
    \else%
      \setlength{\unitlength}{\unitlength * \real{\svgscale}}%
    \fi%
  \else%
    \setlength{\unitlength}{\svgwidth}%
  \fi%
  \global\let\svgwidth\undefined%
  \global\let\svgscale\undefined%
  \makeatother%
  \begin{picture}(1,0.49872378)%
    \lineheight{1}%
    \setlength\tabcolsep{0pt}%
    \put(0,0){\includegraphics[width=\unitlength,page=1]{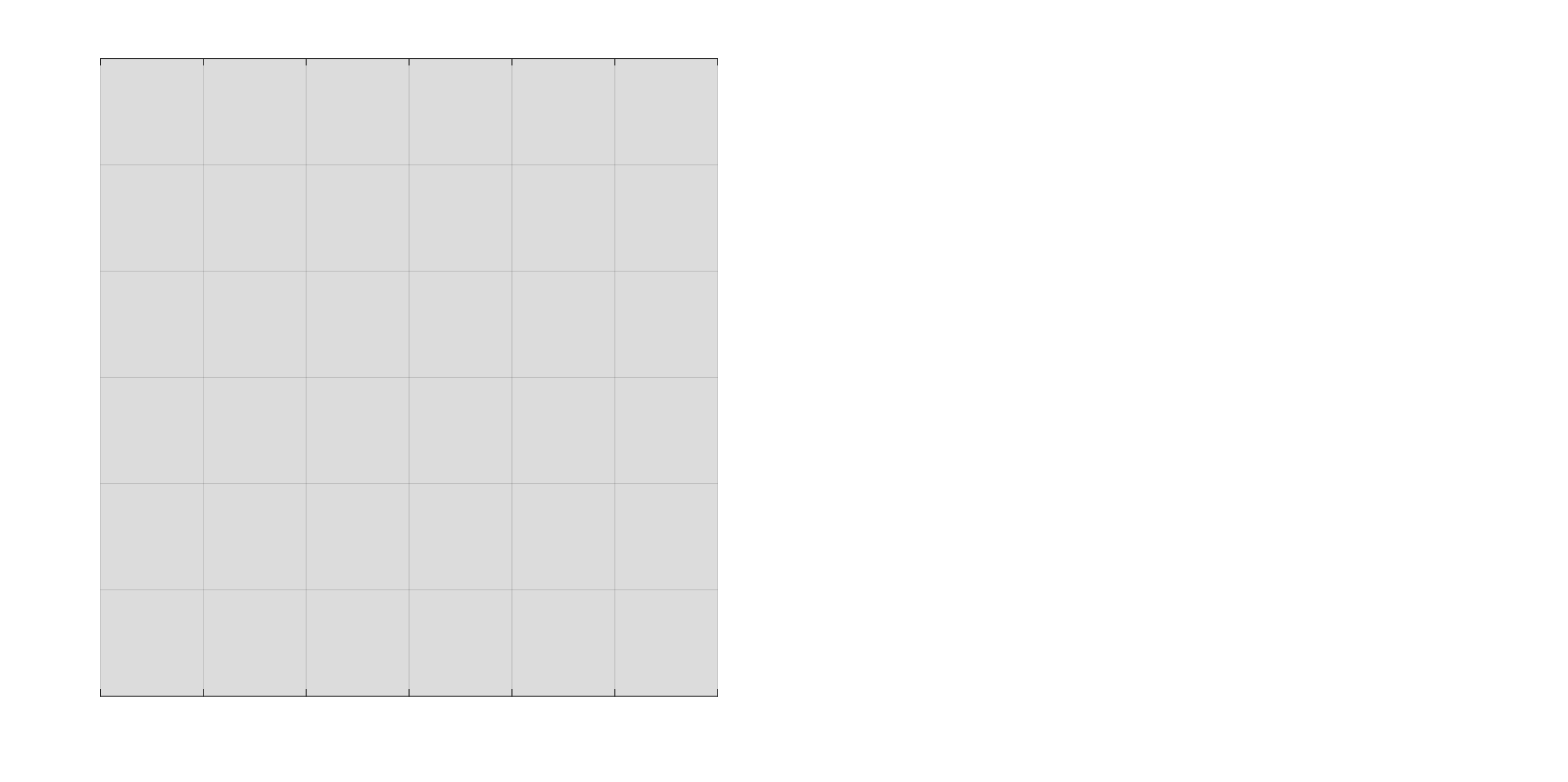}}%
    \put(0.05971486,0.03386169){\makebox(0,0)[lt]{\lineheight{1.25}\smash{\begin{tabular}[t]{l}0\end{tabular}}}}%
    \put(0.11707116,0.03386169){\makebox(0,0)[lt]{\lineheight{1.25}\smash{\begin{tabular}[t]{l}150\end{tabular}}}}%
    \put(0.18270886,0.03386169){\makebox(0,0)[lt]{\lineheight{1.25}\smash{\begin{tabular}[t]{l}300\end{tabular}}}}%
    \put(0.24834655,0.03386169){\makebox(0,0)[lt]{\lineheight{1.25}\smash{\begin{tabular}[t]{l}450\end{tabular}}}}%
    \put(0.31398424,0.03386169){\makebox(0,0)[lt]{\lineheight{1.25}\smash{\begin{tabular}[t]{l}600\end{tabular}}}}%
    \put(0.37962194,0.03386169){\makebox(0,0)[lt]{\lineheight{1.25}\smash{\begin{tabular}[t]{l}750\end{tabular}}}}%
    \put(0.44525963,0.03386169){\makebox(0,0)[lt]{\lineheight{1.25}\smash{\begin{tabular}[t]{l}900\end{tabular}}}}%
    \put(0.20785993,0.01239141){\makebox(0,0)[lt]{\lineheight{1.25}\smash{\begin{tabular}[t]{l}Temperature u\end{tabular}}}}%
    \put(0,0){\includegraphics[width=\unitlength,page=2]{materialparameter.pdf}}%
    \put(0.03591353,0.04876819){\makebox(0,0)[lt]{\lineheight{1.25}\smash{\begin{tabular}[t]{l}30\end{tabular}}}}%
    \put(0.03591353,0.11655291){\makebox(0,0)[lt]{\lineheight{1.25}\smash{\begin{tabular}[t]{l}35\end{tabular}}}}%
    \put(0.03591353,0.18433763){\makebox(0,0)[lt]{\lineheight{1.25}\smash{\begin{tabular}[t]{l}40\end{tabular}}}}%
    \put(0.03591353,0.25212235){\makebox(0,0)[lt]{\lineheight{1.25}\smash{\begin{tabular}[t]{l}45\end{tabular}}}}%
    \put(0.03591353,0.31990707){\makebox(0,0)[lt]{\lineheight{1.25}\smash{\begin{tabular}[t]{l}50\end{tabular}}}}%
    \put(0.03591353,0.38769179){\makebox(0,0)[lt]{\lineheight{1.25}\smash{\begin{tabular}[t]{l}55\end{tabular}}}}%
    \put(0.03591353,0.45547651){\makebox(0,0)[lt]{\lineheight{1.25}\smash{\begin{tabular}[t]{l}60\end{tabular}}}}%
    \put(0.01677432,0.23126552){\rotatebox{90}{\makebox(0,0)[lt]{\lineheight{1.25}\smash{\begin{tabular}[t]{l}k\end{tabular}}}}}%
    \put(0.02474899,0.24598799){\rotatebox{90}{\makebox(0,0)[lt]{\lineheight{1.25}\smash{\begin{tabular}[t]{l}sim\end{tabular}}}}}%
    \put(0.01677432,0.27972699){\rotatebox{90}{\makebox(0,0)[lt]{\lineheight{1.25}\smash{\begin{tabular}[t]{l} (u)\end{tabular}}}}}%
    \put(0,0){\includegraphics[width=\unitlength,page=3]{materialparameter.pdf}}%
    \put(0.57806861,0.03386169){\makebox(0,0)[lt]{\lineheight{1.25}\smash{\begin{tabular}[t]{l}0\end{tabular}}}}%
    \put(0.63552711,0.03386169){\makebox(0,0)[lt]{\lineheight{1.25}\smash{\begin{tabular}[t]{l}150\end{tabular}}}}%
    \put(0.70126712,0.03386169){\makebox(0,0)[lt]{\lineheight{1.25}\smash{\begin{tabular}[t]{l}300\end{tabular}}}}%
    \put(0.76700701,0.03386169){\makebox(0,0)[lt]{\lineheight{1.25}\smash{\begin{tabular}[t]{l}450\end{tabular}}}}%
    \put(0.8327469,0.03386169){\makebox(0,0)[lt]{\lineheight{1.25}\smash{\begin{tabular}[t]{l}600\end{tabular}}}}%
    \put(0.89848692,0.03386169){\makebox(0,0)[lt]{\lineheight{1.25}\smash{\begin{tabular}[t]{l}750\end{tabular}}}}%
    \put(0.96422681,0.03386169){\makebox(0,0)[lt]{\lineheight{1.25}\smash{\begin{tabular}[t]{l}900\end{tabular}}}}%
    \put(0.72652046,0.01239141){\makebox(0,0)[lt]{\lineheight{1.25}\smash{\begin{tabular}[t]{l}Temperature u\end{tabular}}}}%
    \put(0,0){\includegraphics[width=\unitlength,page=4]{materialparameter.pdf}}%
    \put(0.54997322,0.04876819){\makebox(0,0)[lt]{\lineheight{1.25}\smash{\begin{tabular}[t]{l}3.5\end{tabular}}}}%
    \put(0.56224195,0.10686939){\makebox(0,0)[lt]{\lineheight{1.25}\smash{\begin{tabular}[t]{l}4\end{tabular}}}}%
    \put(0.54997322,0.16497059){\makebox(0,0)[lt]{\lineheight{1.25}\smash{\begin{tabular}[t]{l}4.5\end{tabular}}}}%
    \put(0.56224195,0.22307179){\makebox(0,0)[lt]{\lineheight{1.25}\smash{\begin{tabular}[t]{l}5\end{tabular}}}}%
    \put(0.54997322,0.28117292){\makebox(0,0)[lt]{\lineheight{1.25}\smash{\begin{tabular}[t]{l}5.5\end{tabular}}}}%
    \put(0.56224195,0.33927412){\makebox(0,0)[lt]{\lineheight{1.25}\smash{\begin{tabular}[t]{l}6\end{tabular}}}}%
    \put(0.54997322,0.39737532){\makebox(0,0)[lt]{\lineheight{1.25}\smash{\begin{tabular}[t]{l}6.5\end{tabular}}}}%
    \put(0.56224195,0.45547651){\makebox(0,0)[lt]{\lineheight{1.25}\smash{\begin{tabular}[t]{l}7\end{tabular}}}}%
    \put(0.53390119,0.22942521){\rotatebox{90}{\makebox(0,0)[lt]{\lineheight{1.25}\smash{\begin{tabular}[t]{l}C\end{tabular}}}}}%
    \put(0.54187586,0.2478283){\rotatebox{90}{\makebox(0,0)[lt]{\lineheight{1.25}\smash{\begin{tabular}[t]{l}sim\end{tabular}}}}}%
    \put(0.53390119,0.2815673){\rotatebox{90}{\makebox(0,0)[lt]{\lineheight{1.25}\smash{\begin{tabular}[t]{l} (u)\end{tabular}}}}}%
    \put(0.59401795,0.46498478){\makebox(0,0)[lt]{\lineheight{1.25}\smash{\begin{tabular}[t]{l}10\end{tabular}}}}%
    \put(0.62162259,0.47234602){\makebox(0,0)[lt]{\lineheight{1.25}\smash{\begin{tabular}[t]{l}6\end{tabular}}}}%
    \put(0,0){\includegraphics[width=\unitlength,page=5]{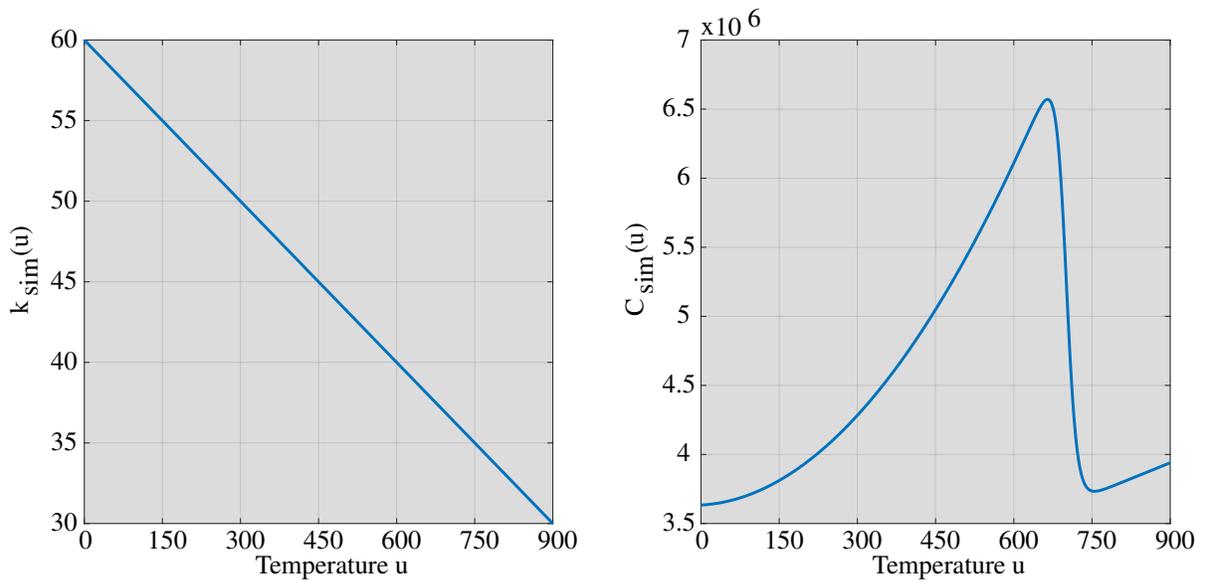}}%
    \put(0.58236266,0.46621165){\color[rgb]{0,0,0}\makebox(0,0)[lt]{\lineheight{1.25}\smash{\begin{tabular}[t]{l}x\end{tabular}}}}%
  \end{picture}%
\endgroup%

\captionof{figure}{Plots of the simulated material parameters $k_{sim}$, $C_{sim} \in \mathcal{C}^1(U)$}
\label{fig:ksimCsim}
\end{center}

Inserting these functions into \eqref{pdeM_Anfang} we prescribe the thermal conduction behaviour of the material. To simulate $M=3$ experiments we specify three different triplets of boundary and initial temperatures, i.e.
\begin{align}\label{triplets}
    \{u^{b,1}(t),u^{t,1}(t),u_0^1(z)\},\ \{u^{b,2}(t),u^{t,2}(t),u_0^2(z)\} \text{ and } \{u^{b,3}(t),u^{t,3}(t),u_0^3(z)\}.
\end{align}

Solving the associated IBVPs \eqref{pdeM_Anfang}-\eqref{pdeM_Ende} numerically, we get the temperature solution matrices $u_i$ for $i=1,\dots,3$, where we extract the three core temperatures, respectively.
To model a plausible measuring process, we add uniformly distributed noise to the data representing an appropriate maximum measurement precision error of $\pm 0.5^{\circ}C$. Thus, our final data will consist of the three noisy core measurements $u^{c,1}$, $u^{c,2}$ and $u^{c,3}$, reflecting the heat conduction behaviour of the material parameters $k_{sim},C_{sim} \in \mathcal{C}^1(U)$. Note, that the components of the data vectors only lie in a subset $\tilde{U}\in U=[0,900]$, which we call the \textit{observed temperature range interval}.

The parameter vector of function values is subsequently determined as solution of the least squares problem

\begin{align} \label{koptCoptFuerSim}
     \left(\underline{k}_{opt},\underline{C}_{opt}\right)^T=\argminA_{\underline{p}:=(\underline{k},\underline{C})^T\in \mathbb{R}^{2n}_+}\sum\limits_{i=1}^{M=3}\left\|QF_i(\underline{p})-u^{c,i}\right\|^2_2
\end{align}

corresponding to the fixed partition $\pi_n$ of $U$. Here, $F_i$ maps a given parameter vector to the solution of \eqref{pdeM_Anfang}-\eqref{pdeM_Ende}, where $C(u)$ and $k(u)$ are the associated PCHIP interpolants. The boundary and inital temperatures are given in \eqref{triplets}. The implementations of $F_1,\ F_2$ and $F_3$ are based on the marching scheme \eqref{marchingScheme}.

For the minimization of the objective funtional in \eqref{koptCoptFuerSim} we use an iterative solver. In particular, we choose a \textit{trust-region-reflective} algorithm of the subroutine \textit{lsqnonlin} provided by the Matlab Optimization Toolbox, see \cite{matlab}. 

 To emphasize that we do not need any a priori information about the functional forms of the material parameters, we choose constant initial guesses, i.e. we have 
 \begin{align}
      k_0(u) &= 45, \\
    C_0(u) &= 4.5*10^6
 \end{align}
 for all $u\in U$. Specifically, this means that we have $\mathbb{R}_+^{2n}\ni\underline{p_0}=(\underline{k_0},\underline{C_0})^T$, where
 \begin{align}
     \underline{k_0}=(\underbrace{45,\dots,45}_{n \text{ times}})^T \text{ and }
     \underline{C_0}=10^6*(\underbrace{4.5,\dots,4.5}_{n \text{ times}})^T
 \end{align}
 are the function values of the PCHIP interpolants $k_0(u)$ and $C_0(u)$ to the partition $\pi_n$ of $U$.
 
 Implementing the forward operators $F_i$ and the observation operator $Q$, the task of finding \eqref{koptCoptFuerSim} consists of a minimization in $\mathbb{R}_+^{2n}$. The iterative solver tries to find a parameter vector $\left(\underline{k}_{opt},\underline{C}_{opt}\right)^T$ by varying the function values and therefore the PCHIP interpolants, like in the example of Figure \ref{fig:pchip}.

The final results, i.e. the calculated PCHIP interpolants $k_{opt}(u)$ and $C_{opt}(u)$ to the initial guesses $k_0(u)$ and $C_0(u)$, are shown in Figure \ref{fig:optResult}.

The small oscillations in $k_{opt}(u)$ and $C_{opt}(u)$ can be neglected and are most likely caused by the solver, i.e. the inherent stopping criteria of the minimizer.
Nevertheless, the form of the functions are retained, but there are gaps between $k_{sim}$ and $k_{opt}$ as well as between $C_{sim}$ and $C_{opt}$. This observation reinforces the fact that the solution of the inverse problem is not unique because the operators $F_i$ are not injective, cf. \eqref{FnichtInjektiv}. As pointed out in subsection \ref{subsection:vierzwei} we can only expect to get a solution pair from a class of solutions having the same quotient. As a result, we can only compare the thermal diffusivities
\begin{align}
    \lambda_{sim}=\frac{k_{sim}}{C_{sim}} \text{\ \ \ and \ \ \   } \lambda_{opt}=\frac{k_{opt}}{C_{opt}},
\end{align}
which yields an outcome, that is illustrated in Figure \ref{fig:quotient} and shows an excellent performance of the method.

\begin{center}
\def\svgwidth{450pt}
\begingroup%
  \makeatletter%
  \providecommand\color[2][]{%
    \errmessage{(Inkscape) Color is used for the text in Inkscape, but the package 'color.sty' is not loaded}%
    \renewcommand\color[2][]{}%
  }%
  \providecommand\transparent[1]{%
    \errmessage{(Inkscape) Transparency is used (non-zero) for the text in Inkscape, but the package 'transparent.sty' is not loaded}%
    \renewcommand\transparent[1]{}%
  }%
  \providecommand\rotatebox[2]{#2}%
  \newcommand*\fsize{\dimexpr\f@size pt\relax}%
  \newcommand*\lineheight[1]{\fontsize{\fsize}{#1\fsize}\selectfont}%
  \ifx\svgwidth\undefined%
    \setlength{\unitlength}{1189.28566007bp}%
    \ifx\svgscale\undefined%
      \relax%
    \else%
      \setlength{\unitlength}{\unitlength * \real{\svgscale}}%
    \fi%
  \else%
    \setlength{\unitlength}{\svgwidth}%
  \fi%
  \global\let\svgwidth\undefined%
  \global\let\svgscale\undefined%
  \makeatother%
  \begin{picture}(1,0.53225227)%
    \lineheight{1}%
    \setlength\tabcolsep{0pt}%
    \put(0,0){\includegraphics[width=\unitlength,page=1]{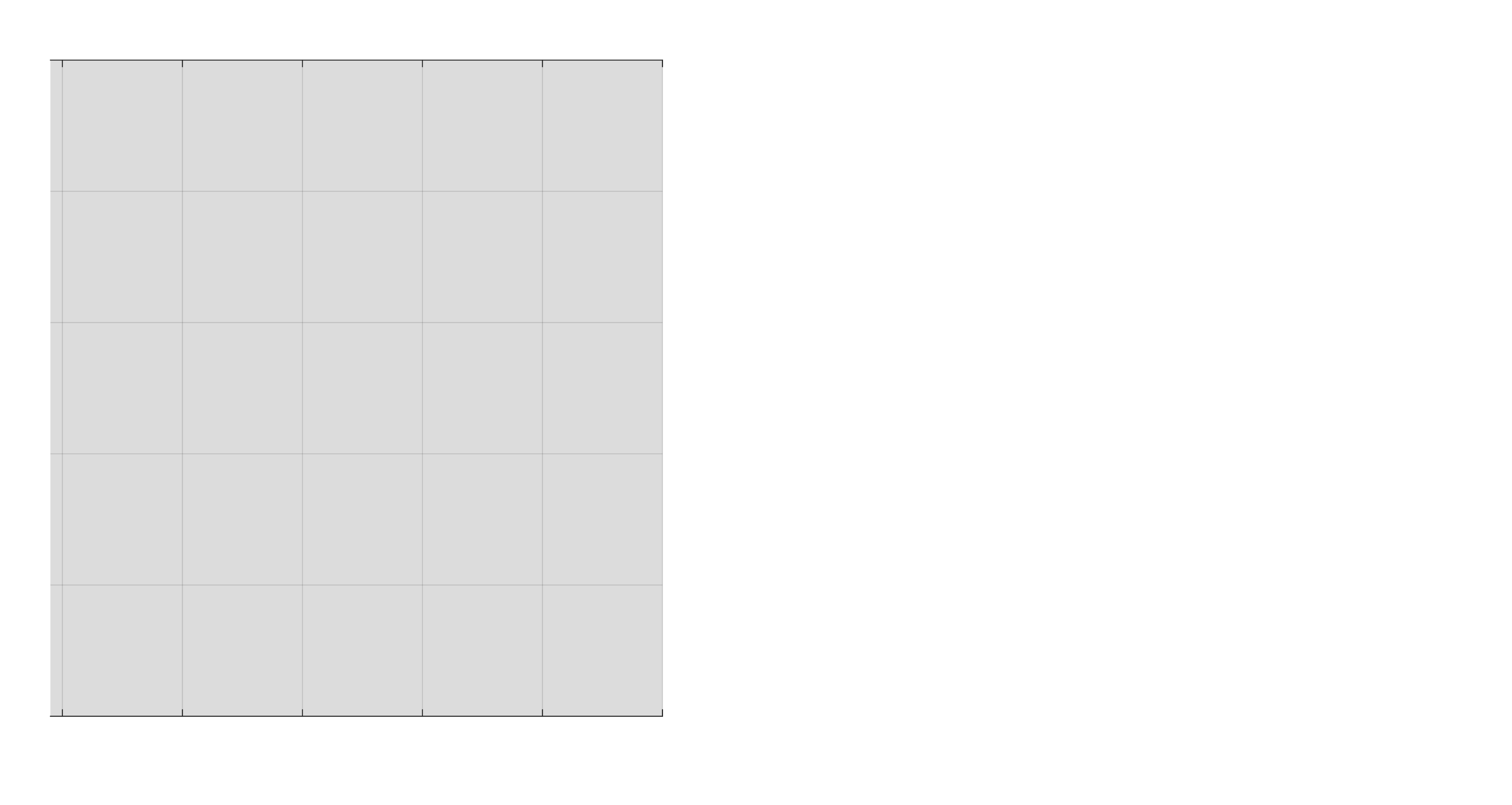}}%
    \put(0.02834393,0.03733333){\makebox(0,0)[lt]{\lineheight{1.25}\smash{\begin{tabular}[t]{l}300\end{tabular}}}}%
    \put(0.10772923,0.03733333){\makebox(0,0)[lt]{\lineheight{1.25}\smash{\begin{tabular}[t]{l}400\end{tabular}}}}%
    \put(0.18711453,0.03733333){\makebox(0,0)[lt]{\lineheight{1.25}\smash{\begin{tabular}[t]{l}500\end{tabular}}}}%
    \put(0.26649977,0.03733333){\makebox(0,0)[lt]{\lineheight{1.25}\smash{\begin{tabular}[t]{l}600\end{tabular}}}}%
    \put(0.34588507,0.03733333){\makebox(0,0)[lt]{\lineheight{1.25}\smash{\begin{tabular}[t]{l}700\end{tabular}}}}%
    \put(0.42527031,0.03733333){\makebox(0,0)[lt]{\lineheight{1.25}\smash{\begin{tabular}[t]{l}800\end{tabular}}}}%
    \put(0.18121643,0.01526126){\makebox(0,0)[lt]{\lineheight{1.25}\smash{\begin{tabular}[t]{l}Temperature u\end{tabular}}}}%
    \put(0,0){\includegraphics[width=\unitlength,page=2]{koptCopt.pdf}}%
    \put(0.00192795,0.05265766){\makebox(0,0)[lt]{\lineheight{1.25}\smash{\begin{tabular}[t]{l}30\end{tabular}}}}%
    \put(0.00192795,0.13943244){\makebox(0,0)[lt]{\lineheight{1.25}\smash{\begin{tabular}[t]{l}35\end{tabular}}}}%
    \put(0.00192795,0.22620721){\makebox(0,0)[lt]{\lineheight{1.25}\smash{\begin{tabular}[t]{l}40\end{tabular}}}}%
    \put(0.00192795,0.31298199){\makebox(0,0)[lt]{\lineheight{1.25}\smash{\begin{tabular}[t]{l}45\end{tabular}}}}%
    \put(0.00192795,0.39975677){\makebox(0,0)[lt]{\lineheight{1.25}\smash{\begin{tabular}[t]{l}50\end{tabular}}}}%
    \put(0.00192795,0.48653155){\makebox(0,0)[lt]{\lineheight{1.25}\smash{\begin{tabular}[t]{l}55\end{tabular}}}}%
    \put(0,0){\includegraphics[width=\unitlength,page=3]{koptCopt.pdf}}%
    \put(0.22882886,0.46603605){\makebox(0,0)[lt]{\lineheight{1.25}\smash{\begin{tabular}[t]{l}k\end{tabular}}}}%
    \put(0.2420721,0.45972974){\makebox(0,0)[lt]{\lineheight{1.25}\smash{\begin{tabular}[t]{l}sim\end{tabular}}}}%
    \put(0,0){\includegraphics[width=\unitlength,page=4]{koptCopt.pdf}}%
    \put(0.31522526,0.46603605){\makebox(0,0)[lt]{\lineheight{1.25}\smash{\begin{tabular}[t]{l}k\end{tabular}}}}%
    \put(0.3284685,0.45972974){\makebox(0,0)[lt]{\lineheight{1.25}\smash{\begin{tabular}[t]{l}0\end{tabular}}}}%
    \put(0,0){\includegraphics[width=\unitlength,page=5]{koptCopt.pdf}}%
    \put(0.39090094,0.46603605){\makebox(0,0)[lt]{\lineheight{1.25}\smash{\begin{tabular}[t]{l}k\end{tabular}}}}%
    \put(0.40414418,0.45972974){\makebox(0,0)[lt]{\lineheight{1.25}\smash{\begin{tabular}[t]{l}opt\end{tabular}}}}%
    \put(0,0){\includegraphics[width=\unitlength,page=6]{koptCopt.pdf}}%
    \put(0.56123919,0.03733333){\makebox(0,0)[lt]{\lineheight{1.25}\smash{\begin{tabular}[t]{l}300\end{tabular}}}}%
    \put(0.64074816,0.03733333){\makebox(0,0)[lt]{\lineheight{1.25}\smash{\begin{tabular}[t]{l}400\end{tabular}}}}%
    \put(0.720257,0.03733333){\makebox(0,0)[lt]{\lineheight{1.25}\smash{\begin{tabular}[t]{l}500\end{tabular}}}}%
    \put(0.79976597,0.03733333){\makebox(0,0)[lt]{\lineheight{1.25}\smash{\begin{tabular}[t]{l}600\end{tabular}}}}%
    \put(0.87927487,0.03733333){\makebox(0,0)[lt]{\lineheight{1.25}\smash{\begin{tabular}[t]{l}700\end{tabular}}}}%
    \put(0.95878384,0.03733333){\makebox(0,0)[lt]{\lineheight{1.25}\smash{\begin{tabular}[t]{l}800\end{tabular}}}}%
    \put(0.71441471,0.01526126){\makebox(0,0)[lt]{\lineheight{1.25}\smash{\begin{tabular}[t]{l}Temperature u\end{tabular}}}}%
    \put(0,0){\includegraphics[width=\unitlength,page=7]{koptCopt.pdf}}%
    \put(0.53039644,0.05265766){\makebox(0,0)[lt]{\lineheight{1.25}\smash{\begin{tabular}[t]{l}3.5\end{tabular}}}}%
    \put(0.54300905,0.11463963){\makebox(0,0)[lt]{\lineheight{1.25}\smash{\begin{tabular}[t]{l}4\end{tabular}}}}%
    \put(0.53039644,0.17662161){\makebox(0,0)[lt]{\lineheight{1.25}\smash{\begin{tabular}[t]{l}4.5\end{tabular}}}}%
    \put(0.54300905,0.23860358){\makebox(0,0)[lt]{\lineheight{1.25}\smash{\begin{tabular}[t]{l}5\end{tabular}}}}%
    \put(0.53039644,0.30058562){\makebox(0,0)[lt]{\lineheight{1.25}\smash{\begin{tabular}[t]{l}5.5\end{tabular}}}}%
    \put(0.54300905,0.3625676){\makebox(0,0)[lt]{\lineheight{1.25}\smash{\begin{tabular}[t]{l}6\end{tabular}}}}%
    \put(0.53039644,0.42454957){\makebox(0,0)[lt]{\lineheight{1.25}\smash{\begin{tabular}[t]{l}6.5\end{tabular}}}}%
    \put(0.54300905,0.48653155){\makebox(0,0)[lt]{\lineheight{1.25}\smash{\begin{tabular}[t]{l}7\end{tabular}}}}%
    \put(0.58450455,0.49756758){\makebox(0,0)[lt]{\lineheight{1.25}\smash{\begin{tabular}[t]{l}10\end{tabular}}}}%
    \put(0.60909914,0.50639641){\makebox(0,0)[lt]{\lineheight{1.25}\smash{\begin{tabular}[t]{l}6\end{tabular}}}}%
    \put(0,0){\includegraphics[width=\unitlength,page=8]{koptCopt.pdf}}%
    \put(0.61414419,0.46792794){\makebox(0,0)[lt]{\lineheight{1.25}\smash{\begin{tabular}[t]{l}C\end{tabular}}}}%
    \put(0.63054058,0.46162164){\makebox(0,0)[lt]{\lineheight{1.25}\smash{\begin{tabular}[t]{l}sim\end{tabular}}}}%
    \put(0,0){\includegraphics[width=\unitlength,page=9]{koptCopt.pdf}}%
    \put(0.70369374,0.46792794){\makebox(0,0)[lt]{\lineheight{1.25}\smash{\begin{tabular}[t]{l}C\end{tabular}}}}%
    \put(0.72009014,0.46162164){\makebox(0,0)[lt]{\lineheight{1.25}\smash{\begin{tabular}[t]{l}0\end{tabular}}}}%
    \put(0,0){\includegraphics[width=\unitlength,page=10]{koptCopt.pdf}}%
    \put(0.78252257,0.46792794){\makebox(0,0)[lt]{\lineheight{1.25}\smash{\begin{tabular}[t]{l}C\end{tabular}}}}%
    \put(0.79891897,0.46162164){\makebox(0,0)[lt]{\lineheight{1.25}\smash{\begin{tabular}[t]{l}opt\end{tabular}}}}%
    \put(0,0){\includegraphics[width=\unitlength,page=11]{koptCopt.pdf}}%
    \put(0.56747752,0.49630633){\color[rgb]{0,0,0}\makebox(0,0)[lt]{\lineheight{1.25}\smash{\begin{tabular}[t]{l}x\end{tabular}}}}%
  \end{picture}%
\endgroup%

\captionof{figure}{Comparison of exact (simulated), initial and calculated (optimized) material parameters} in the observed temperature range interval $\tilde{U}$
\label{fig:optResult}
\end{center}

\begin{center}
\def\svgwidth{400pt}
\begingroup%
  \makeatletter%
  \providecommand\color[2][]{%
    \errmessage{(Inkscape) Color is used for the text in Inkscape, but the package 'color.sty' is not loaded}%
    \renewcommand\color[2][]{}%
  }%
  \providecommand\transparent[1]{%
    \errmessage{(Inkscape) Transparency is used (non-zero) for the text in Inkscape, but the package 'transparent.sty' is not loaded}%
    \renewcommand\transparent[1]{}%
  }%
  \providecommand\rotatebox[2]{#2}%
  \newcommand*\fsize{\dimexpr\f@size pt\relax}%
  \newcommand*\lineheight[1]{\fontsize{\fsize}{#1\fsize}\selectfont}%
  \ifx\svgwidth\undefined%
    \setlength{\unitlength}{1439.99998186bp}%
    \ifx\svgscale\undefined%
      \relax%
    \else%
      \setlength{\unitlength}{\unitlength * \real{\svgscale}}%
    \fi%
  \else%
    \setlength{\unitlength}{\svgwidth}%
  \fi%
  \global\let\svgwidth\undefined%
  \global\let\svgscale\undefined%
  \makeatother%
  \begin{picture}(1,0.4234375)%
    \lineheight{1}%
    \setlength\tabcolsep{0pt}%
    \put(0,0){\includegraphics[width=\unitlength,page=1]{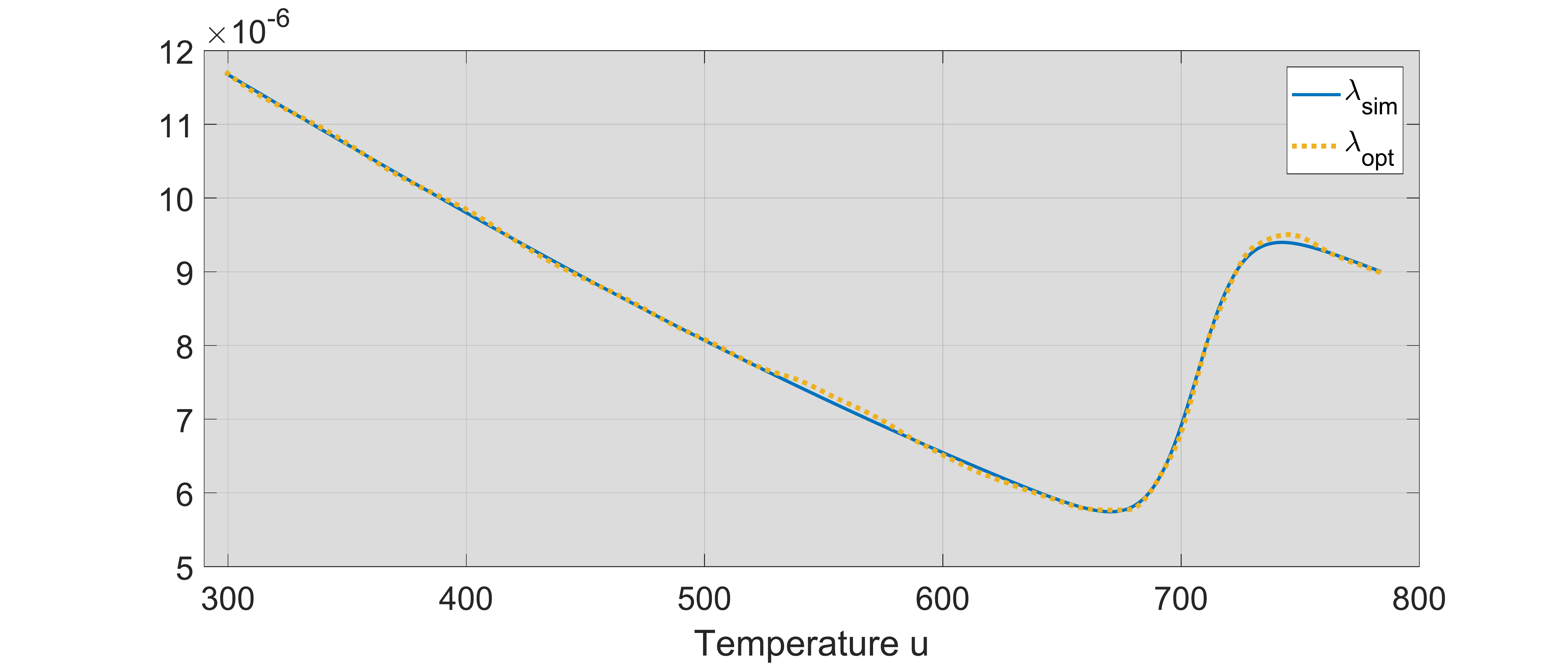}}%
  \end{picture}%
\endgroup%

\captionof{figure}{Comparison of exact (simulated) and calculated (optimized) thermal diffusivities}
\label{fig:quotient}
\end{center}

\section{Conclusion}
We were able to numerically determine the temperature dependent material parameters $k(u)$ and $C(u)$ (up to some canonical ambiguity) in a 1D nonlinear heat equation describing the heat conduction over the thickness of a cooled heavy plate. With the general interpolation procedure described in this paper, we can break down the identification process of functions from an infinite-dimensional space $\mathcal{C}^1(U)$ to a minimization in $\mathbb{R}_+^{2n}$.  Although we added noise to the data, we can show that the conduction behaviour of the material,  in terms of the thermal diffusivity, can be reconstructed almost perfectly. Besides that, we did not need any a priori information about the functional shape of the material parameters. Our investigations show that the method used in this paper is applicable to real data in the production of TMCP steel plates. There, the heat conduction determination provides physically plausible and valuable results in  characterizing the underlying material, which is the first step in understanding and controlling the Accelerated Cooling (ACC) process. Future work is focusd on modeling the heat fluxes on the surfaces depending on the surface temperature itself (Leidenfrost effect) and the control variables of the cooling device, i.e. the water load in every cooling zone and the feed rate of the heavy plate.
 

\bibliographystyle{elsarticle-num}
\bibliography{references_arXiv}

\begin{thebibliography}{10}
\expandafter\ifx\csname url\endcsname\relax
  \def\url#1{\texttt{#1}}\fi
\expandafter\ifx\csname urlprefix\endcsname\relax\def\urlprefix{URL }\fi
\expandafter\ifx\csname href\endcsname\relax
  \def\href#1#2{#2} \def\path#1{#1}\fi

\bibitem{egger}
H.~Egger, J.-F. Pietschmann, M.~Schlottbom, Identification of nonlinear heat
  conduction laws, Journal of Inverse and Ill-posed Problems 23~(5) (2015)
  429--437.

\bibitem{cannon}
J.~Cannon, Determination of the unknown coefficient $k (u)$ in the equation
  $\nabla\cdot (k (u) \nabla u)= 0$ from overspecified boundary data, Journal
  of Mathematical Analysis and Applications 18~(1) (1967) 112--114.

\bibitem{cannonDuchateau}
J.~Cannon, P.~Duchateau, Determining unknown coefficients in a nonlinear heat
  conduction problem, SIAM Journal on Applied Mathematics 24~(3) (1973)
  298--314.

\bibitem{rincon1}
M.~Rincon, J.~L{\'\i}maco, I.~S. Liu, Existence and uniqueness of solutions of
  a nonlinear heat equation, Trends in Applied and Computational Mathematics
  6~(2) (2005) 273--284.

\bibitem{rincon2}
M.~Teixeira, M.~Rincon, I.-S. Liu, Numerical analysis of quenching--heat
  conduction in metallic materials, Applied Mathematical Modelling 33~(5)
  (2009) 2464--2473.

\bibitem{benyuZou}
G.~Ben-yu, J.~Zou, An augmented lagrangian method for parameter identifications
  in parabolic systems, Journal of mathematical analysis and applications
  263~(1) (2001) 49--68.

\bibitem{engl}
H.~W. Engl, J.~Zou, A new approach to convergence rate analysis of tikhonov
  regularization for parameter identification in heat conduction, Inverse
  Problems 16~(6) (2000) 1907.

\bibitem{hussein}
M.~Hussein, D.~Lesnic, M.~Ivanchov, Simultaneous determination of
  time-dependent coefficients in the heat equation, Computers \& Mathematics
  with Applications 67~(5) (2014) 1065--1091.

\bibitem{kunisch}
K.~Kunisch, G.~Peichl, Estimation of a temporally and spatially varying
  diffusion coefficient in a parabolic system by an augmented lagrangian
  technique, Numerische Mathematik 59~(1) (1991) 473--509.

\bibitem{hankeScherzer}
M.~Hanke, O.~Scherzer, Error analysis of an equation error method for the
  identification of the diffusion coefficient in a quasi-linear parabolic
  differential equation, SIAM Journal on Applied Mathematics 59~(3) (1998)
  1012--1027.

\bibitem{cuiGaoZhang}
M.~Cui, X.~Gao, J.~Zhang, A new approach for the estimation of
  temperature-dependent thermal properties by solving transient inverse heat
  conduction problems, International Journal of Thermal Sciences 58 (2012)
  113--119.

\bibitem{huangYan}
C.-H. Huang, Y.~Jan-Yuan, An inverse problem in simultaneously measuring
  temperature-dependent thermal conductivity and heat capacity, International
  Journal of Heat and Mass Transfer 38~(18) (1995) 3433--3441.

\bibitem{roubicek}
T.~Roub{\'\i}{\v{c}}ek, Nonlinear partial differential equations with
  applications, Vol. 153, Springer Science \& Business Media, 2013.

\bibitem{louis}
A.~K. Louis, Inverse and ill-posed problems, Springer-Verlag, 2013.

\bibitem{schuster}
T.~Schuster, B.~Kaltenbacher, B.~Hofmann, K.~S. Kazimierski, Regularization
  methods in Banach spaces, Vol.~10, Walter de Gruyter, 2012.

\bibitem{fritschCarlson}
F.~N. Fritsch, R.~E. Carlson, Monotone piecewise cubic interpolation, SIAM
  Journal on Numerical Analysis 17~(2) (1980) 238--246.

\bibitem{matlab}
Matlab optimization toolbox 8.0, the MathWorks, Natick, MA, USA.

\end{thebibliography}


\end{document}